\documentclass[11pt,a4paper,reqno]{amsart} 

\usepackage[dvipsnames]{xcolor}
\usepackage{tikz}
\usepackage{amsmath,bm,bbm,amsthm, amssymb}
\usepackage{fullpage}
\usepackage{color}
\usepackage{caption}
\usepackage{subcaption}

\usepackage[hyperindex,breaklinks]{hyperref}
\hypersetup{ colorlinks=true, linkcolor=blue, citecolor=blue,filecolor=blue, urlcolor=blue}

\definecolor{wiasblue}   {cmyk}{1.0, 0.60, 0, 0}
\definecolor{mlugreen}{RGB}{172,6,52}

\newcommand\bq[1]{{\bf  #1}}

\def\ni{\noindent}
\def\Z{\mathbb Z}
\def\E{\mathbb E}
\def\P{\mathbb P}
\def\K{\mathbb K}

\def\R{\mathbb R}
\def\mc{\mathcal}
\def\ms{\mathsf}

\def\la{\lambda}
\def\La{\Lambda}
\def\a{\alpha}
\def\aa{\boldsymbol a}
\def\bb{\boldsymbol b}
\def\dd{\boldsymbol d}

\def\s{\sigma}
\def\su{\subseteq}
\def\bs{\boldsymbol}
\def\e{\varepsilon}
\def\eb{\varepsilon_{\ms B}}

\def\g{\gamma}
\def\b{\beta}
\def\de{\delta}

\def\om{\omega}

\def\es{\varnothing}
\def\one{\mathbbmss{1}}
\def\ccv{C_{\ms{CV}}}

\def\De{\Delta}

\def\ff{\infty}

\def\vp{\varphi}

\def\d{{\rm d}}
\def\k{\kappa}

\def\GG{\mc G}

\def\PP{X}

\def\f{\frac}

\def\off{[0,\ff)}

\def\r{\rho}

\def\im{\item}
\def\sm{\setminus}

\def\bef{\begin{figure}}
\def\enf{\end{figure}}
\def\bep{\begin{proof}}
\def\enp{\end{proof}}
\def\bepr{\begin{proposition}}
\def\enpr{\end{proposition}}
\def\bec{\begin{corollary}}
\def\enc{\end{corollary}}
\def\bea{\begin{align}}
\newcommand\eea{\end{align}}
\def\beas{\begin{align*}}
\def\eeas{\end{align*}}
\def\bet{\begin{theorem}}
\def\ent{\end{theorem}}
\def\bee{\begin{example}}
\def\ene{\end{example}}

\def\bede{\begin{definition}}
\def\ende{\end{definition}}
\def\ber{\begin{remark}}
\def\enr{\end{remark}}
\def\beca{\begin{cases}}
\def\enca{\end{cases}}
\def\bel{\begin{lemma}}
\def\enl{\end{lemma}}
\def\been{\begin{enumerate}}
\def\enen{\end{enumerate}}

\def\SS{\mc S}
\def\dim{\ms{dim}}
\def\beit{\begin{itemize}}
\def\enit{\end{itemize}}
\def\befr{\begin{frame}}
\def\enfr{\end{frame}}
\def\ti{\times}

\def\L{\mathbb L}
\def\Var{\ms{Var}}
\def\Cov{\ms{Cov}}

\def\ba{|\,}

\def\bi{\big}

\def\diam{\ms{diam}}
\def\ker{\ms{ker}}
\def\Im{\ms{Im}}
\def\wt{\widetilde}

\def\Cech{\ms{ Cech}}

\renewcommand\le{\leqslant}
\renewcommand\ge{\geqslant}

\def\dist{\ms{dist}}

\def\becb{\begin{tcolorbox}[colback=Dandelion!20]}
\def\encb{\end{tcolorbox}}

\def\co{\colon}
\def\bn{\b_n^{\bb, \dd}}
\def\bnp{\b_n^{\bb', \dd'}}

\def\bbn{\bar\b_n^{\bb, \dd}}
\def\bqn{\b_{q, n}^{\bb, \dd}}

\def\bqp{\b_q^{\bb', \dd'}}
\def\bq{\b_q^{\bb, \dd}}
\def\bbq{\bar\b_q^{\bb, \dd}}

\def\bab{\bar\b}

\def\Ebbb{E_{\ms b,3}}

\def\Ebe{E_{\ms b,m}}

\def\Edbb{E_{\ms d,3}}
\def\pa{\partial}

\def\Edd{E_{\ms d,m}}

\def\Var{\ms{Var}}
\def\kn{n^d}

\def\dnz{\De_{i, n}}

\def\bnz{\b_{i, n}}
\def\dnz{\De_{i, n}}
\def\Mu{\ms{Mult}}
\def\MC{\K}
\def\SSS{\mc S}
\def\rC{r}
\def\oto{\SSS}
\def\ig{\includegraphics}
\def\kb{k_{\ms b}}
\def\kd{k_{\ms d}}

\theoremstyle{plain}
\newtheorem{theorem}{Theorem}[section]
\newtheorem{proposition}[theorem]{Proposition}
\newtheorem{corollary}[theorem]{Corollary}
\newtheorem{lemma}[theorem]{Lemma}

\theoremstyle{definition}
\newtheorem{definition}[theorem]{Definition}

\newtheorem{example}[theorem]{Example}

\theoremstyle{remark}
\newtheorem{remark}[theorem]{Remark}

\keywords{topological data analysis, persistence diagram,  multi-parameter persistence, Goodness-of-fit tests , consistency, asymptotic normality}
\subjclass[2010]{60F05; 60D05; 60G55 ; 55U10}
\date{\today}

\begin{document}

\title{On the consistency and asymptotic normality of multiparameter persistent Betti numbers}

 \author{Magnus B.~Botnan}
        \author{Christian Hirsch}
        \address[Magnus B. Botnan]{Vrije Universiteit, Department of Mathematics, Faculteit der Exacte Wetenschappen, De Boelelaan 1111, 1081 HV Amsterdam}
       \email{M.B.Botnan@vu.nl}
        \address[Christian Hirsch]{University of Groningen, Bernoulli Institute, Nijenborgh 9, 9747 AG Groningen, The Netherlands.}
       \email{c.p.hirsch@rug.nl}

\begin{abstract}
	The persistent Betti numbers are used in topological data analysis to infer the scales at which topological features appear and disappear in the filtration of a topological space. Most commonly by means of the corresponding barcode or persistence diagram. While this approach to data science has been very successful, it suffers from sensitivity to outliers, and it does not allow for additional filtration parameters. Such parameters naturally appear when a cloud of data points comes together with additional measurements taken at the locations of the data. For these reasons, multiparameter persistent homology has recently received significant attention. In particular, the multicover and \v Cech bifiltration have been introduced to overcome the aforementioned shortcomings.

In this work, we establish the strong consistency and asymptotic normality of the multiparameter persistent Betti numbers in growing domains. Our asymptotic results are established for a general framework encompassing both the marked \v Cech bifiltration, as well as the multicover bifiltration constructed on the null model of an independently marked Poisson point process. In a simulation study, we explain how the asymptotic normality can be used to derive tests for the goodness of fit. The statistical power of such tests is illustrated through different alternatives exhibiting more clustering, or more repulsion than the null model.
\end{abstract}

\maketitle
\section{Introduction}

\begin{figure}
\begin{subfigure}{.32\textwidth}
  \centering
  \includegraphics[width=.9\linewidth]{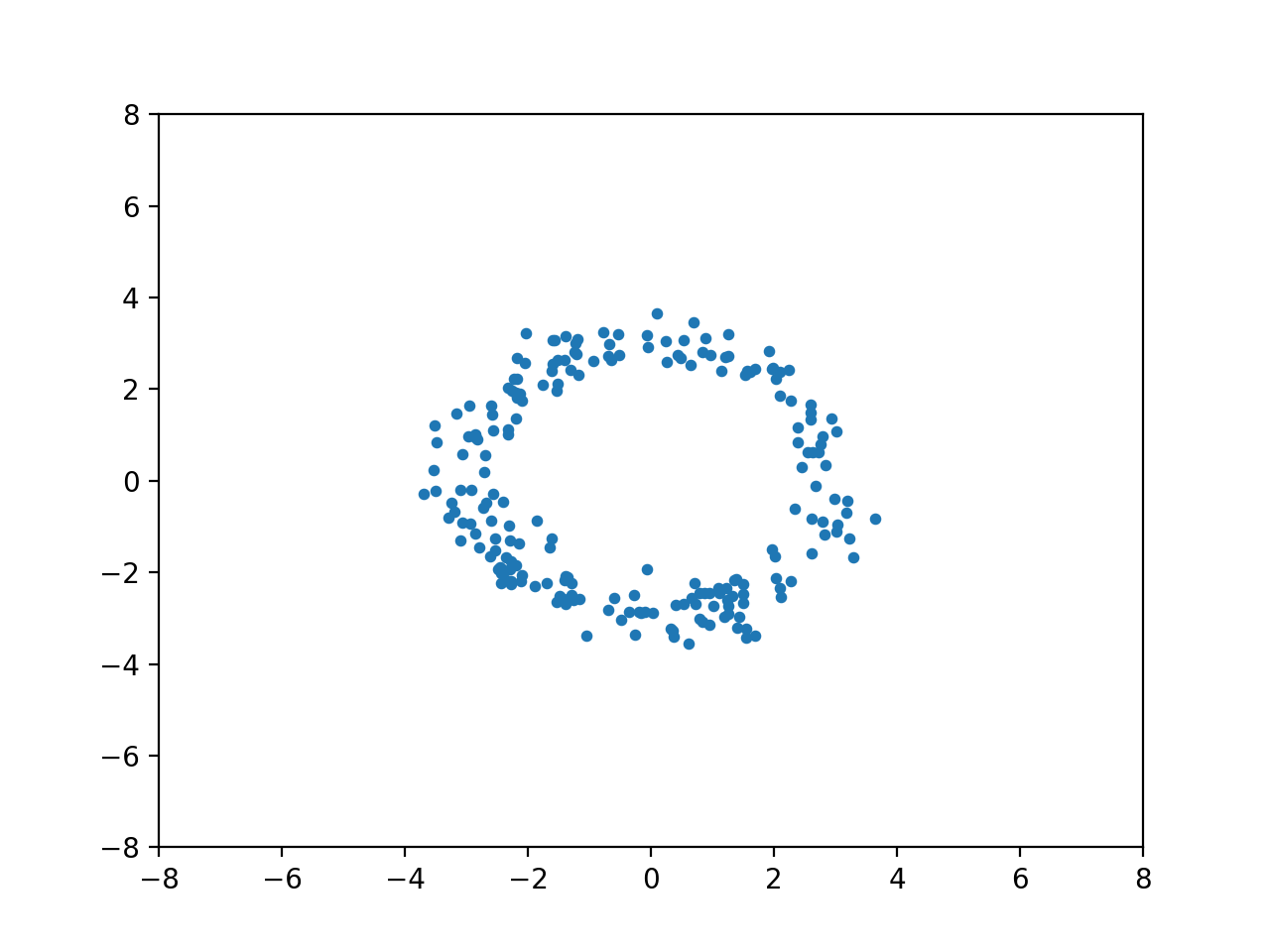}
  \caption{}
  \label{fig:circle}
\end{subfigure}%
\begin{subfigure}{.32\textwidth}
  \centering
  \includegraphics[width=.9\linewidth]{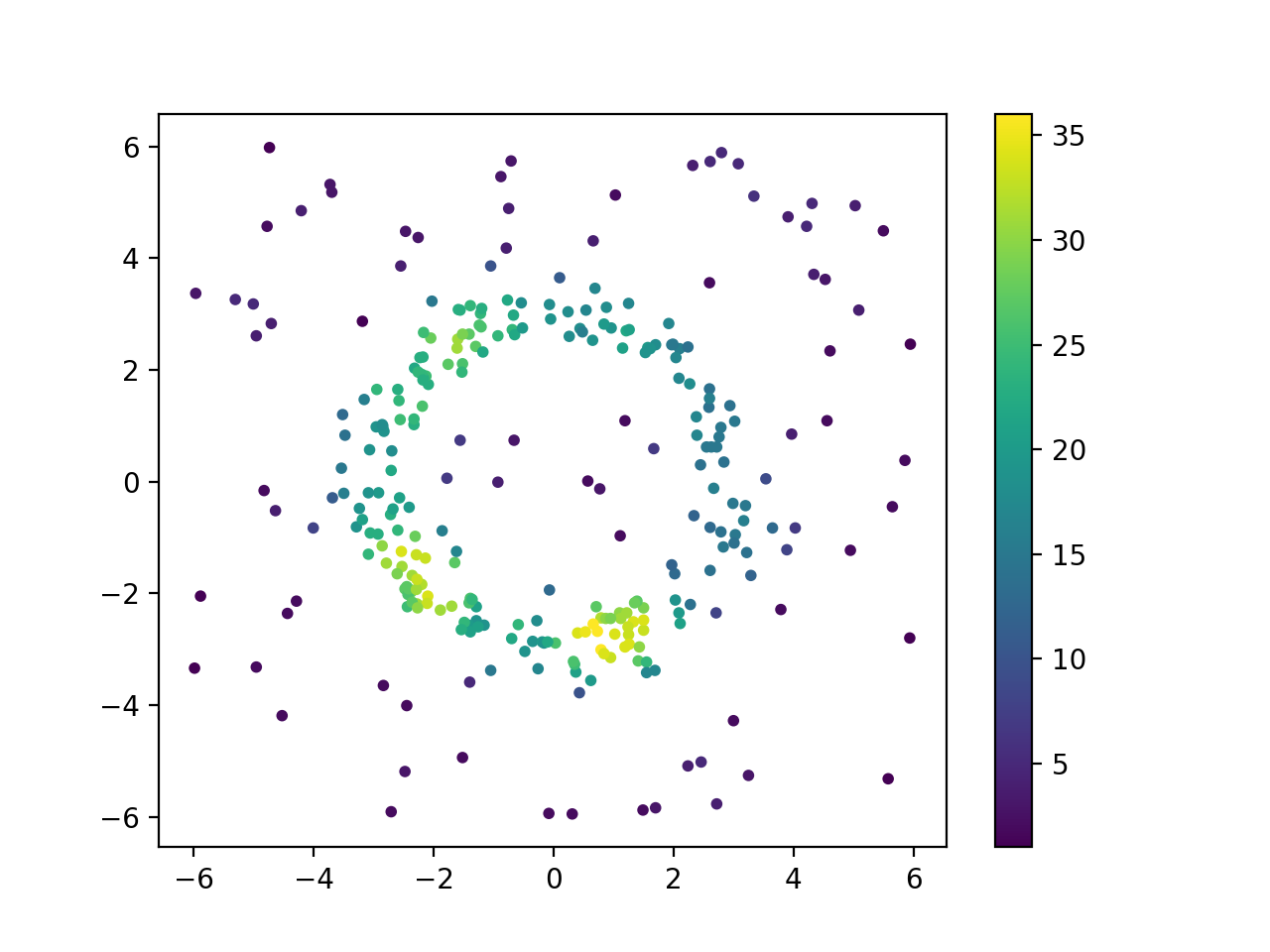}
  \caption{}
  \label{fig:circle-all-points}
\end{subfigure}
\begin{subfigure}{.32\textwidth}
  \centering
  \includegraphics[width=.9\linewidth]{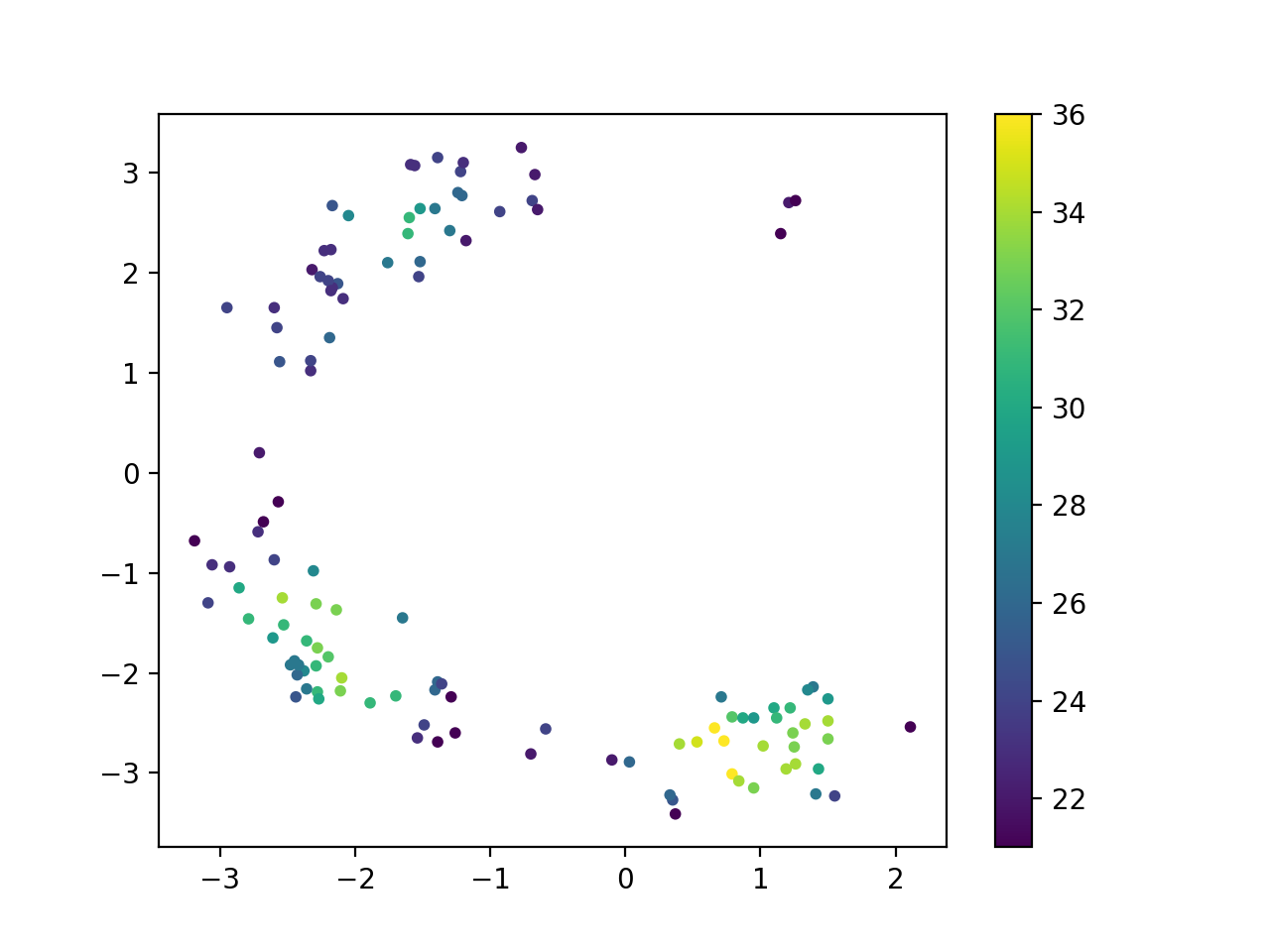}
  \caption{}
  \label{fig:circle-too-few}
\end{subfigure}
\begin{subfigure}{.32\textwidth}
  \centering
  \includegraphics[width=.9\linewidth]{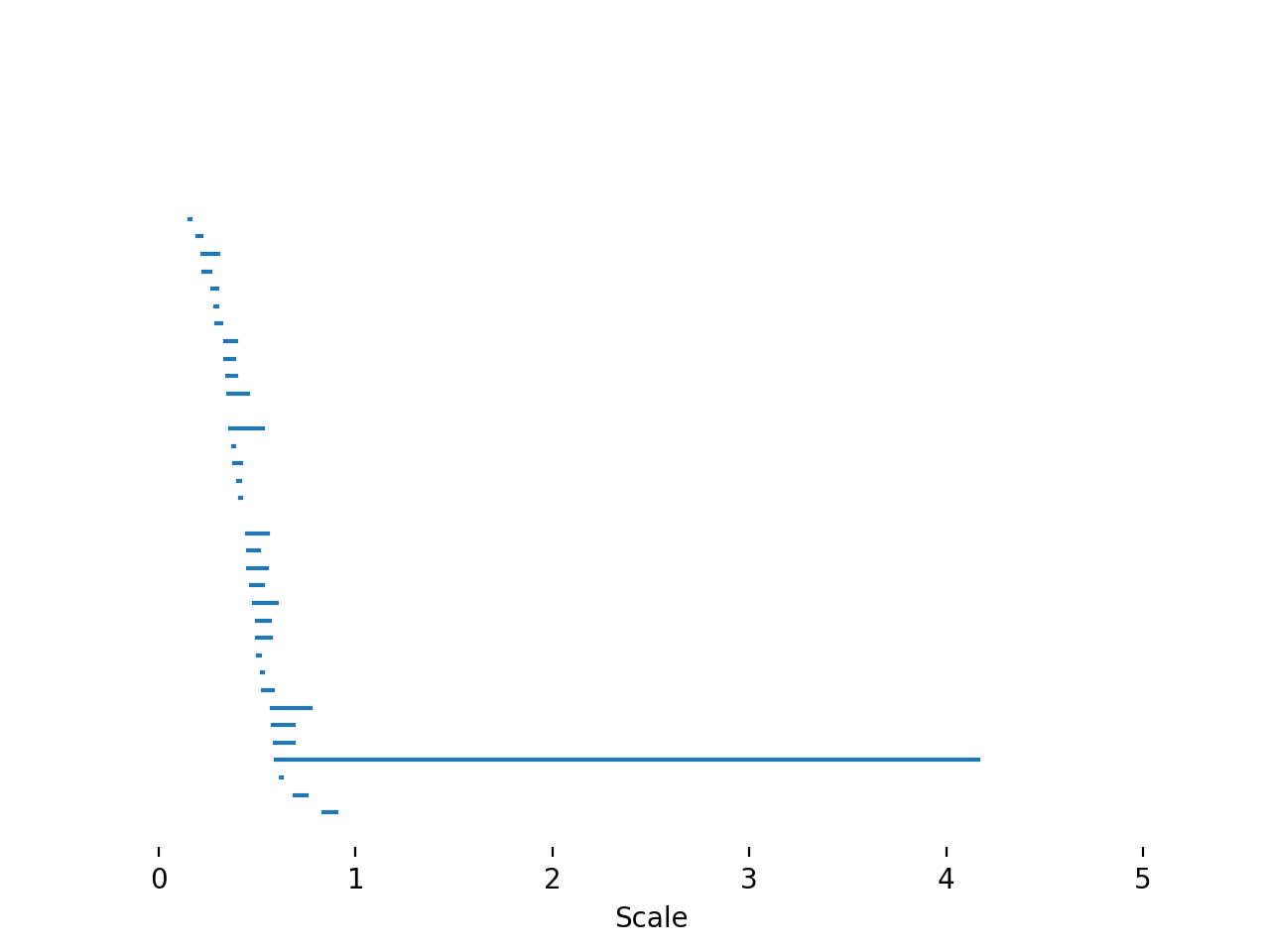}
  \caption{}
  \label{fig:barcode-circle}
\end{subfigure}%
\begin{subfigure}{.32\textwidth}
  \centering
  \includegraphics[width=.9\linewidth]{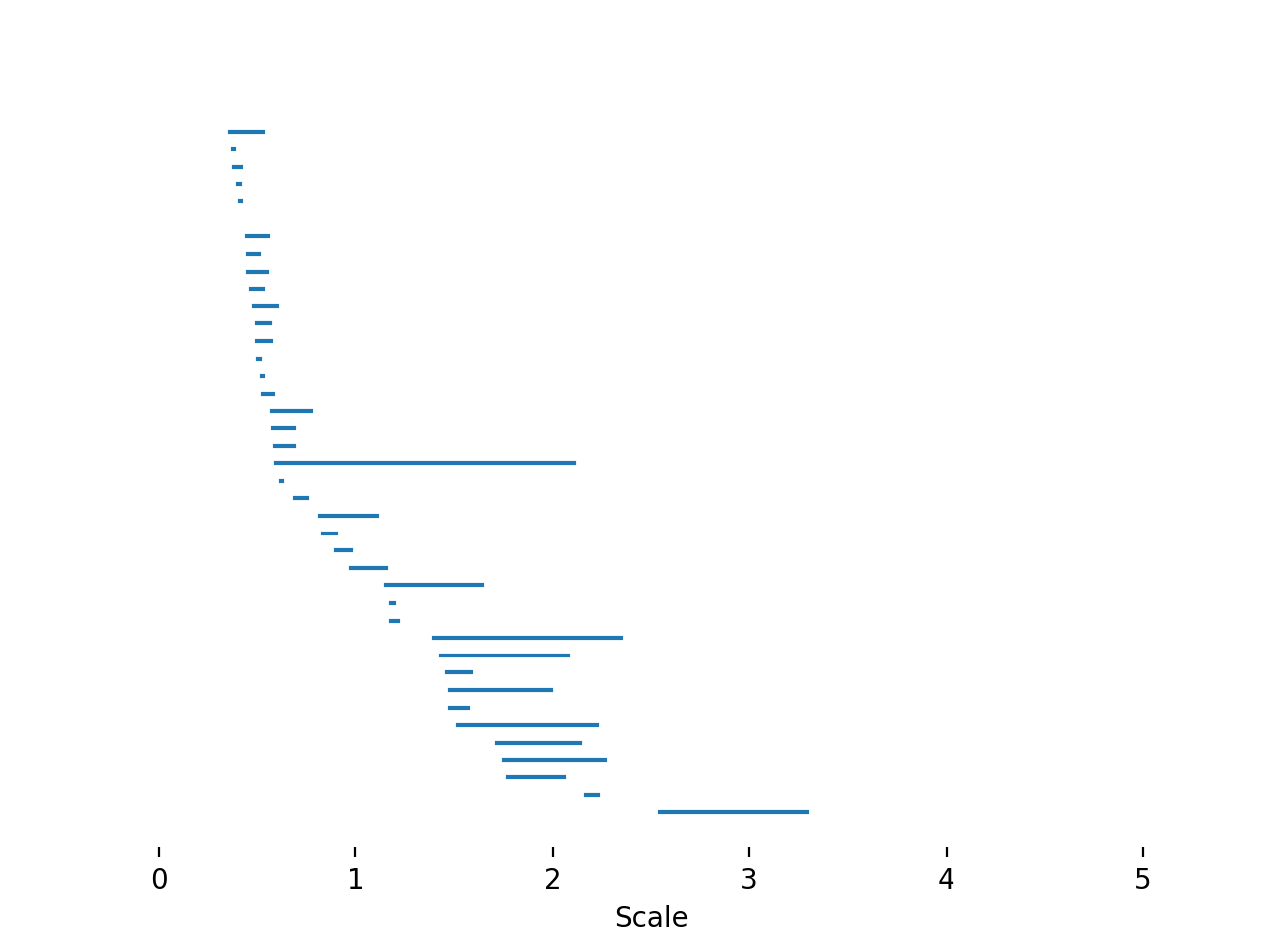}
  \caption{}
  \label{fig:barcode-all}
\end{subfigure}
\begin{subfigure}{.32\textwidth}
  \centering
  \includegraphics[width=.9\linewidth]{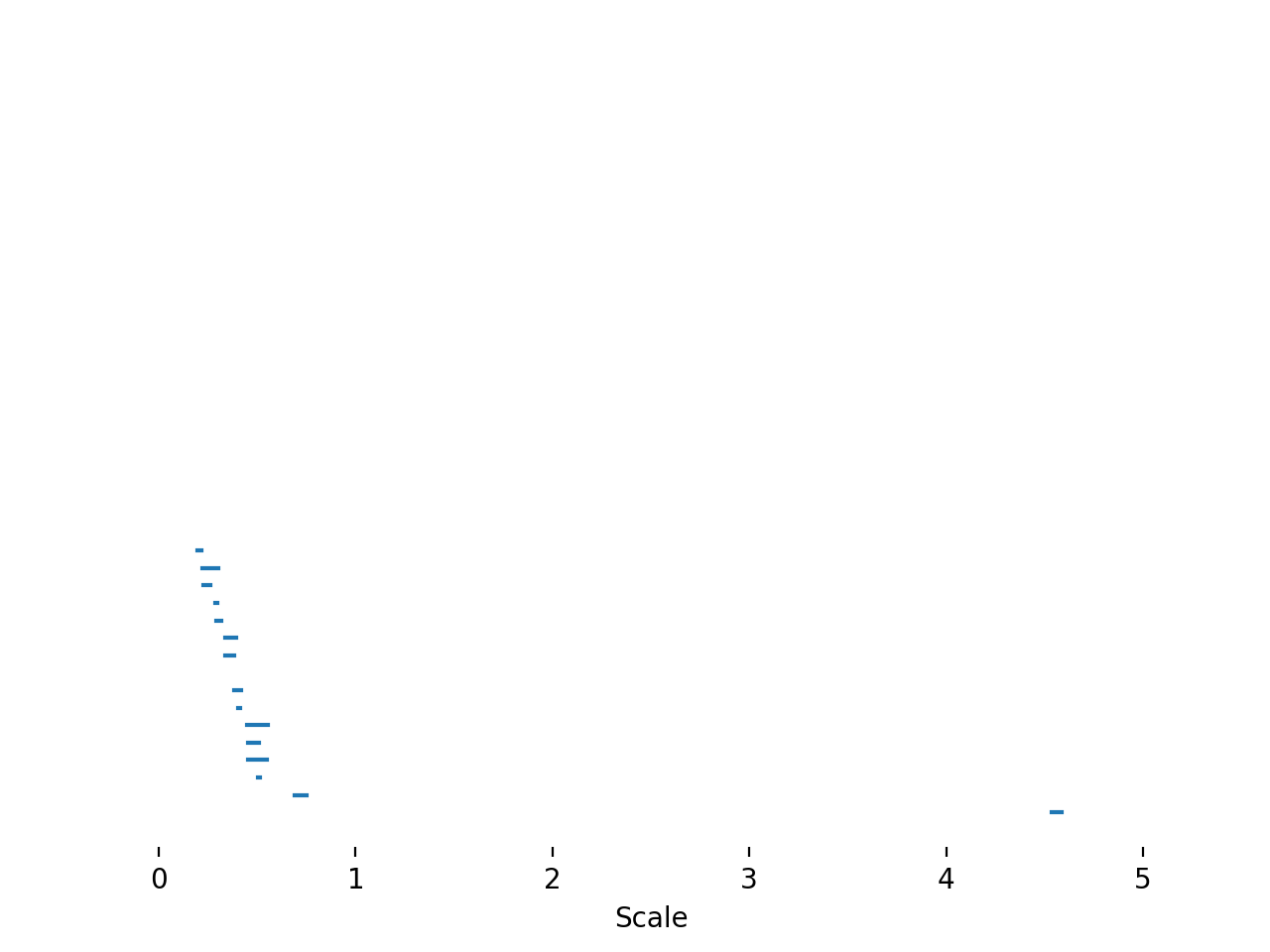}
  \caption{}
  \label{fig:barcode-few}
\end{subfigure}
\caption{\textbf{(a)} A data set with a circular shape. \textbf{(b)} The data from (a) with added noise and colorized by a local density estimate. \textbf{(c)} The data points in (b) with a sufficiently high local density estimate. \textbf{(d)} The barcode of the data in (a). \textbf{(e)} The barcode of the data in (b). \textbf{(f)} The barcode of the data in (c). }
\end{figure}

The goal of topological data analysis (TDA) is to infer the ''shape'' of data by means of topological invariants. One of the most notable such tools, \emph{persistent homology}, outputs a collection of intervals, the \emph{barcode}, that track the evolution of the topological (homological) features along a filtration of a simplicial complex. Short intervals are treated  as ''noise'' and long bars as true topological signals, with the precise notion of what constitutes a long interval being application-dependent.  
For instance, the circular feature of the data set in  Figure \ref{fig:circle} is readily deduced from the associated barcode shown in  \ref{fig:barcode-circle}. Note that it is also customary to visualize the collection of intervals as a \emph{persistence diagram} in which the interval $[a,b)$ in the barcode is plotted as the point $(a,b)\in \R^2$; see Figure \ref{pd_mc_fig} for an example. 

The persistence diagram has now become a powerful tool to analyze complex phenomena in fields as diverse as astronomy, biology, materials science, medicine and neuroscience \cite{wasserman}. This rapid dissemination in the natural sciences has lead to a vigorous research stream aiming to put TDA on a rigorous statistical foundation \cite{bendich,divol,jmva}.

What makes persistent homology particularly appealing to the data-analyst is the fact that it is stable: data sets that are close to each other (e.g. in the Hausdorff distance) also have similar barcodes. At the same time, a drawback of persistent homology is its sensitivity to outliers. To see this, consider the same circle as before but with a few additional points scattered around as in Figure \ref{fig:circle-all-points}. The associated barcode in Figure \ref{fig:barcode-all} no longer suggests an underlying circular structure. A potential way to rectify this problem would be to consider only points above a certain density threshold, but this would in turn be very sensitive to the choice of threshold. See for instance Figure \ref{fig:circle-too-few} and its associated barcode in Figure \ref{fig:barcode-few}. Ideally one would thus have a tool which allows one to deduce topological signatures across both scale and density. That is precisely the promise of \emph{multi-parameter persistent homology} is about.

While constructing multi-parameter persistent homology is straightforward,  the transition from a single to multiple parameters comes at the price of a significant complexity increase of the underlying algebraic objects. This has prompted the introduction of novel invariants \cite{harrington2019stratifying,miller2020homological} but so far much of the work has be centered around the idea of studying the family of barcodes obtained by restricting the indexing set to a set of straight lines \cite{cerri2013betti,rivet}. The information contained in such restrictions is equivalent to the data of the \emph{rank invariant}: the collection of ranks between every pair of comparable points.  Since the rank invariant is one of few efficiently computable invariants \cite{botnan_et_al:LIPIcs:2020:12180,carlsson2009theory}, it offers a natural starting point for the development of a sound statistical foundation for multi-parameter persistent homology. 

%RESULTS
More precisely, for a prototypical model in multi-parameter persistence, we prove a functional law of large numbers (LLN) and a functional central limit theorem (CLT) in large domains, thereby showing the strong consistency and asymptotic normality of a flexible class of test statistics. Our general framework comprises two examples that are of particular interest. 

%MARK CECH
Our first example for bifiltered simplicial complex -- the \emph{marked \v Cech bifiltration} -- concerns locations that are scattered at random in a sampling window according to a Poisson point process, and that are endowed with some independent marking. For instance, we may think of the points as measurement locations with the mark as the quantity that is being measured.  Indeed, in the axis of the point locations, we can consider the standard \v Cech filtration. Additionally, we can use sub-level sets of the mark space to define a second filtration, which effectively leads to a thinning of the measurement locations whose mark does not exceed a certain value.

%MULTICOV
Our second example, the \emph{multicover bifiltration} from \cite{osang} connects to the topic of robustness to outliers mentioned above. Also this filtration builds on locations scattered in a sampling window. The 1-cover corresponds to the ordinary \v Cech complex, which describes the topology of the union of growing balls centered at the locations. In general, in the $k$-cover we extract more refined information which captures the topology of the set covered by at least $k \ge 2$ of the balls. Hence, the corresponding characteristics will not be influenced by the occurrence of isolated outliers located in atypical locations.

%METHOD
The present work is the first step towards a rigorous statistical foundation of multi-parameter persistence. Moreover, also on the methodological level, we introduce novel proof ideas in order to deal with the challenges in the multi-parameter setting in comparison to the established results in the literature \cite{krebs}. The most fundamental difference is that in the multi-parameter setting, there is no general analog of the persistence diagram. Thus, it is no longer possible to interpret the increments of the rank invariant in terms of topological features with birth- and death times, which requires a more detailed analysis of the relevant geometric configurations. Moreover, we need a precise control over the H\"older continuity of the marks in order to establish the moment bounds in the tightness part of the CLT proof.

%ORG
The rest of the manuscript is organized as follows. First, Section \ref{mod_sec} contains a precise definition of the multi-parameter persistence model described above. Then, Section \ref{res_sec} presents the main results of this work, namely a functional strong LLN and a functional CLT for the rank invariant, thus yielding asymptotic normality for a flexible class of test statistics. In a simulation study in Section \ref{sim_sec}, we illustrate how to leverage the functional CLT in order to develop specific goodness-of-fit tests for different point patterns. Finally, in Section \ref{conc_sec}, we summarize the findings and provide an outlook to further research. The detailed proofs for the main results are then given in Sections \ref{prel_sec}.

\section{Model}
\label{mod_sec}
We shall assume that the reader is familiar with simplicial complexes and simplicial homology. We recommend \cite{wasserman} for a highly accessible account of these topics, which is written for an audience with a background in statistics. Henceforth, $\K$ always denotes a \emph{multifiltered simplicial complex}. That is, 
\been
\im $\K = \{K_{\aa}\}_{\aa}$ describes a family of simplicial complexes indexed by $\aa = (a_1, \dots, a_u) \in \SSS := \mc S_1 \times \cdots \times \mc S_u$ for some totally ordered sets $\mc S_1, \dots, \mc S_u$;
\im for every $\aa = (a_1, \dots, a_u), \bb = (b_1, \dots, b_u) \in \SSS$ with $a_i \le b_i$ for every $i \le u$ there is a simplicial map $K_{\aa} \to K_{\bb}$.  Here, a simplicial map is a map from the 0-simplices of $K_{\aa}$ to the 0-simplices of $K_{\bb}$ such that every $q$-simplex in $K_{\aa}$ is mapped to a $q$-simplex in $K_{\bb}$.
\enen

%
%POINT PROCESS
%
In this work, we consider bifiltrations built on geometric point patterns. More precisely, throughout the entire manuscript, we fix some deterministic $T > 0$ and let $ X = \{(X_i, M_i)\}_{i \ge 1}$ denotes a $[0, T]$-marked point process, which is stationary in the sense that the distribution of $\{(X_i + x, M_i)\}_{i \ge 1}$ is the same for every $x \in \R^d$. In applications, the mark could allow for a variety of different interpretations such as the volume of a  particle, the amount of current flowing through a given measurement location, or the total precipitation at a given location. We assume that the intensity $\la := \E\big[\#\{i\co X_i \in [0, 1]^d\}\big]$ is positive and finite. For instance, $\{X_i\}_{i \ge 1}$ can be a homogeneous Poisson point process endowed with iid marks \cite{poisBook}.

The conceptual framework developed in the present article is guided by two fundamental examples of bifiltrations: the \emph{marked \v Cech-bifiltration} and the \emph{multicover bifiltration}  \cite{osang, marked}. In order to render the presentation of the overarching framework more accessible, we discuss these examples first before moving to the general setting. 
Loosely speaking, the \emph{marked \v Cech bifiltration (\v C-bifiltration)} illustrated in Figure \ref{bifilt_fig} combines the \v C-filtration with the sub-level filtration on the marks.
%
%MARKED
%
\bee[Marked \v Cech bifiltration]
\label{mark_exc}
Set
\begin{align}
	\label{bifi_eq}
	\K_{(r_1, r_2), n} := \Cech_{r_1}\big(\{X_i \in [0, n]^d \co M_i \le r_2\}\big).
\end{align}
Here, the \emph{filtration time} in the \v C-filtration of a $q$-simplex $\{x_0, \dots, x_q\} \su \R^d$ is given by 
$$\rC(\{x_0, \dots, x_q\}) := \min\bi\{t > 0\co B_t(x_0) \cap \cdots \cap B_t(x_q) \ne \es\bi\},$$
where $B_t(y) := \{x \in \R^d\co |x - y| \le t\}$ denotes the Euclidean ball with radius $t > 0$ centered at $y \in \R^d$.
The simplicial maps are given by the natural inclusions of simplicial complexes. 

\bef[!h]
\centering
	        \input{bifilt}
			\caption{Marked \v C-bifiltration with uniform marks on $[0, 1]$. Illustration shows the bifiltration at levels $(0.8, 0.7)$ (left), $(0.8, 0.2)$ (center), and $(1.2, 0.2)$ (right). For clarity, only Delaunay triangles are drawn.}
			\label{bifilt_fig}
\enf

%
%MULT MARK
%
Note that it would also be possible to consider multi-variate marks, i.e., $M_i \in [0, T]^p$ for some $p \ge 1$. In order to render the presentation more accessible, we focus on the case $p = 1$, noting that most of the arguments in this work extends to $p \ge 1$ after straightforward modifications.
\ene

%
%MULTICOV
%
\bee[Multicover bifiltration]
\label{cov_exc}
Although highly popular in applications, the \v C-filtration has the drawback of being sensitive to outliers \cite{merigot}. A small number of data points placed in sparse regions in space may influence sensitively the topological properties of the corresponding \v C-complex. One attractive option to mitigate these effects is to work with \emph{multicover filtrations}. Instead of studying the topology of the union of balls centered at data points, one works with the set of points that are covered by a certain number $k \ge 1$ of such balls. Figure \ref{vis_fig} illustrates the multiple coverage for a set of random points in a 2D sampling window.
\bef[!h]
\centering
		\input{mult_vis}
			\caption{Areas of multiple coverage shown by gray shades. Small radius (left) and large radius (right).}
			\label{vis_fig}
\enf

Some work is needed to embed the idea of multicover bifiltrations into the general framework of simplicial complexes. For $k \ge 1$, $r > 0$, the 0-simplices of the complex $\Mu_{(r, k)}(X)$ are given by all unordered $k$-tuples of points in $\{X_i\}_{i \ge 1}$. Then, $\vp_0, \dots, \vp_p \su \{X_i\}_{i \ge 1}$ with $\#\vp_i = k$ form a $p$-simplex at level $r > 0$ if and only if $\vp_0 \cup \cdots \cup \vp_p$ forms a simplex in the standard \v Cech filtration. In symbols, if
\begin{align}
	\label{mc_eq}
\vp_0 \cup \cdots \cup \vp_p \in \Cech_r\big(\{X_i \in [0, n]^d\}\big).
\end{align}

In order to turn $\Mu_{(r, k)}(X)$ into a bifiltration, we reverse the natural ordering on the $k$-component. Hence, in order to define a simplicial map from $\Mu_{(r_1, k_1)}(X)$ to $\Mu_{(r_2, k_2)}(X)$ for $r_1 \le r_2$ and $k_1 \ge k_2$, we need to provide a map sending a subset $\vp_1 \su \{X_i\}_{i \ge 1}$ of size $k_1$ to a subset $\vp_2 \su \{X_i\}_{i \ge 1}$ of size $k_2$. This is achieved by defining $\vp_2$ to consist of the $k_2$ smallest elements of $\vp_1$ according to the lexicographic order. This choice is of course arbitrary but we stress that a different selection rule does not alter the persistence properties of the resulting complex.
\ene

Finally, it is possible to unify Examples \ref{mark_exc} and \ref{cov_exc} under a single umbrella. Loosely speaking, we compute the multicover bifiltration only among those points whose mark is smaller than a certain filtration level.
%
%COMBINED
%
\bee[Combined example]
\label{comb_exc}
Set
\begin{align}
	\label{comb_eq}
\MC_{(r_1, r_2, k), n}(\{(X_i, M_i)\}_{i \ge 1}) := \Mu_{(r_1, k)}\big(\{X_i \in [0, n]^d \co M_i \le r_2\}\big).
\end{align}
Henceforth, we fix some $T', K' \ge 1$ and work on the compact index set  $\SSS := [0, T'] \ti [0, T] \ti \{1, \dots, K'\}$. In order to ease the presentation, we assume from now on that $T = T' = K'$.  In particular, setting $k = 1$ recovers the  {marked \v Cech-bifiltration}, and setting $r_2 = T$ recovers the  {multicover bifiltration}.
\ene

%
%DIAG
%
For single-parameter filtrations, the persistence diagram is a collection of pairs $\{(B_j, D_j)\}_j$ representing the filtration values where features such as connected components, loops, or higher-dimensional cavities appear and disappear. For instance by fixing one of the axes and varying the other one, Figure \ref{pd_mark_fig} provides two persistence diagrams that can be extracted from the marked \v C-bifiltration.
\bef[!h]
\centering
	        \ig[width=.4\textwidth]{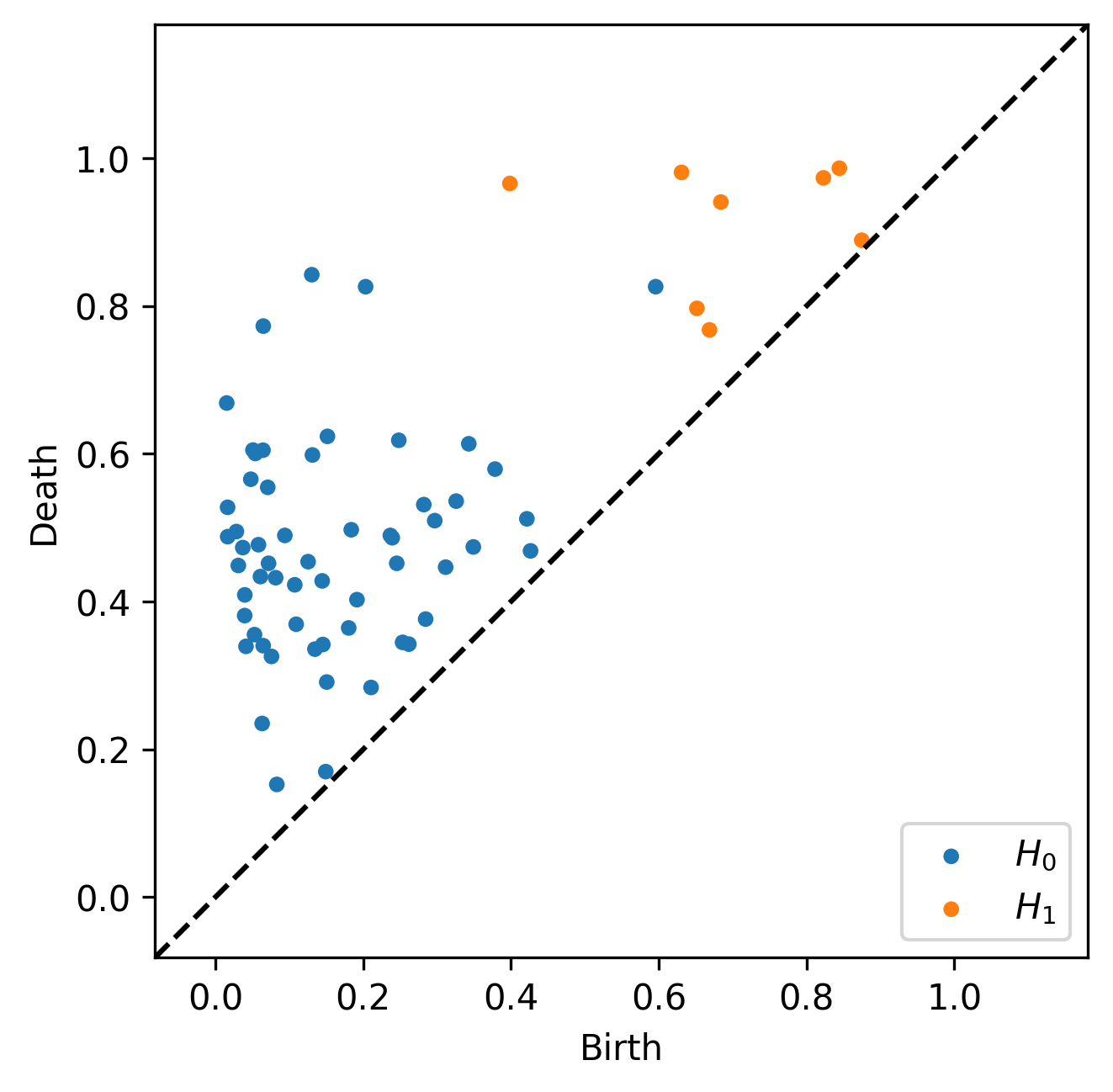}
		        \ig[width=.4\textwidth]{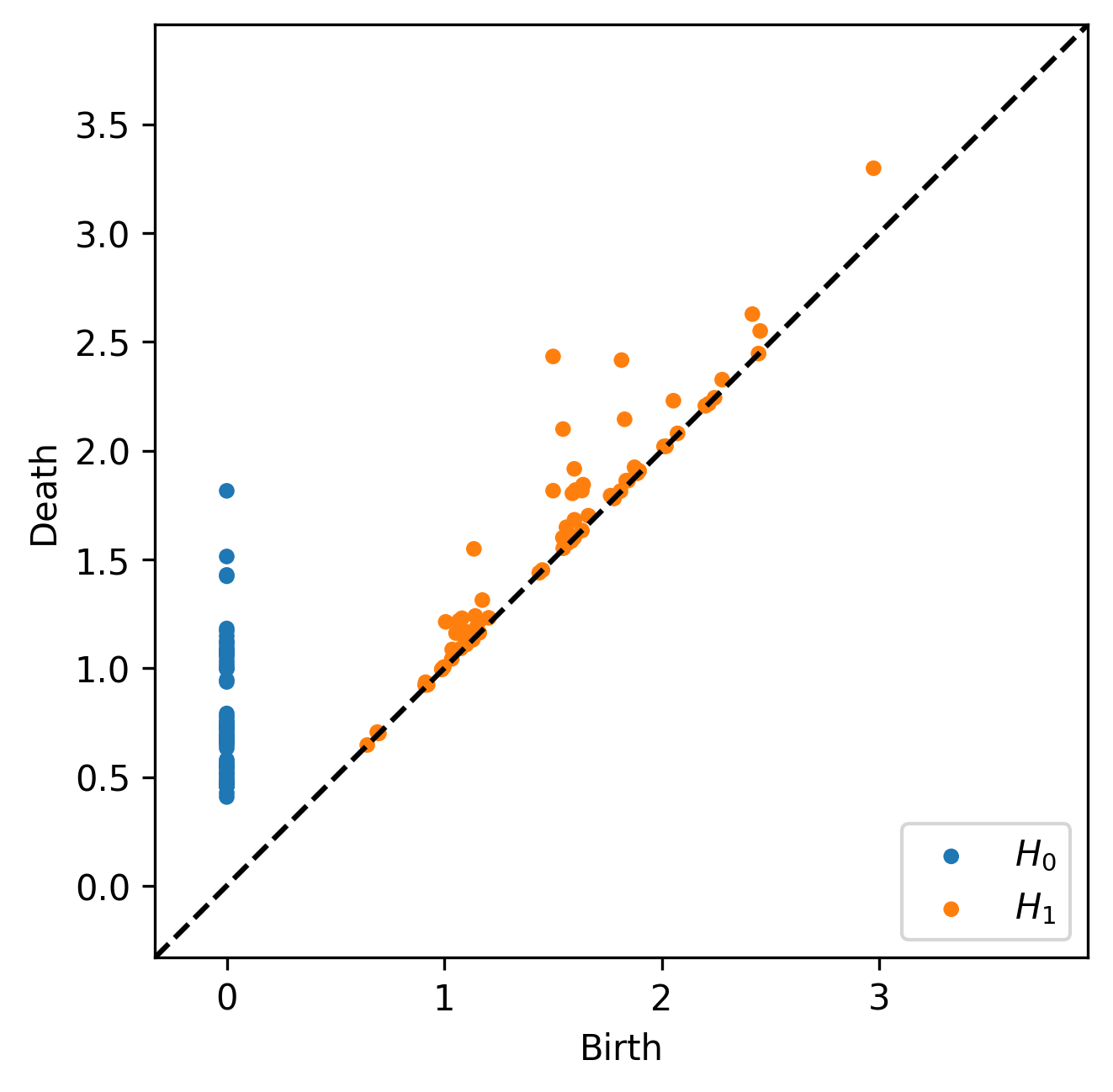}
\caption{Persistence diagrams constructed from the point clouds in Figure \ref{bifilt_fig} with respect to the mark-axis (left) and \v Cech-axis (right).}
\label{pd_mark_fig}
\enf

Similarly, the multicover bifiltration yields persistence diagrams for different fixed values of the covering depth $k \ge 1$. For instance, Figure \ref{pd_mc_fig} illustrates these diagrams when $k = 1$ and $k = 2$. 
\bef[!h]\centering
	        \ig[width=.4\textwidth]{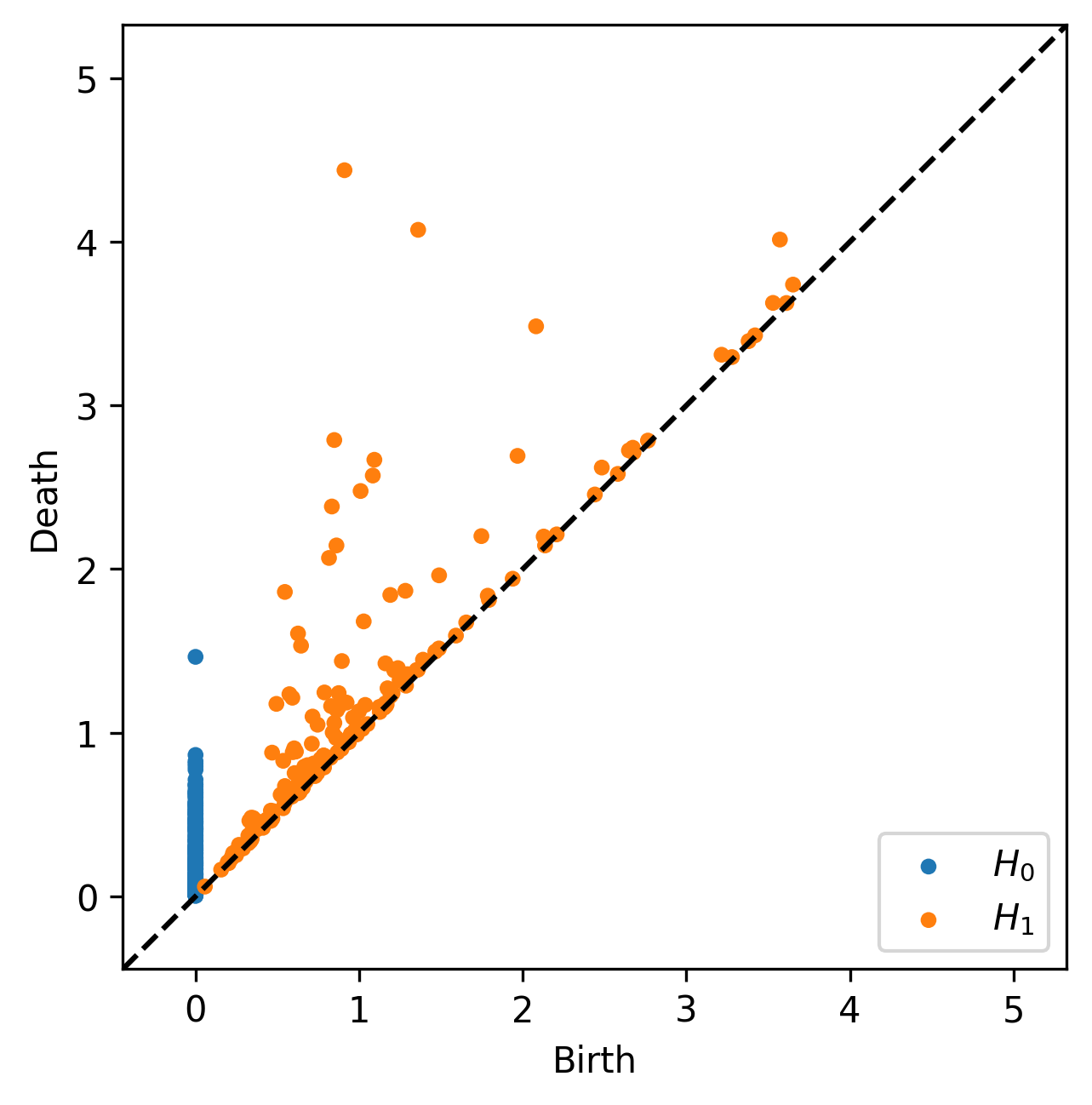}
		        \ig[width=.4\textwidth]{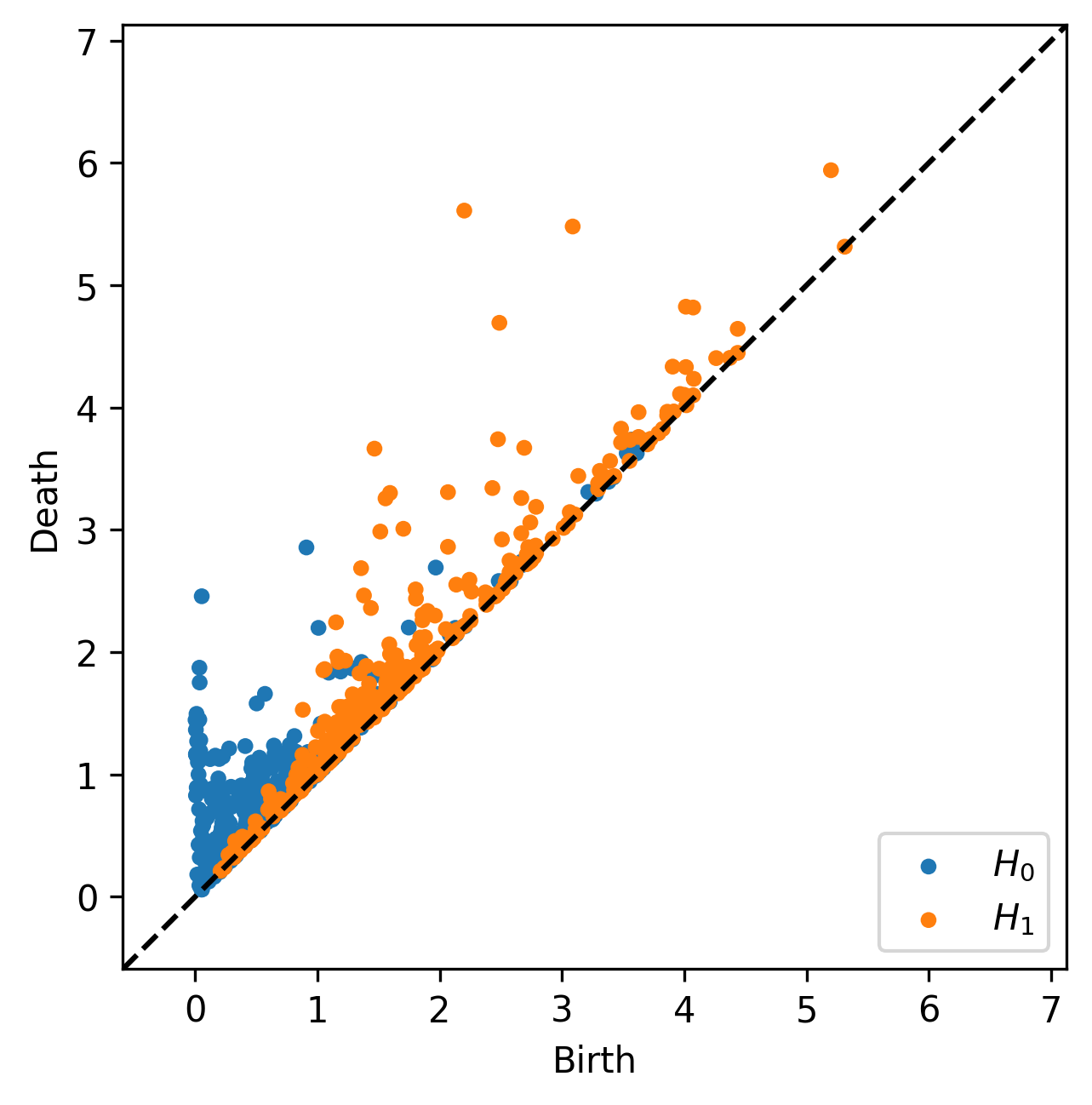}
				\caption{Persistence diagrams constructed from the point clouds in Figure \ref{vis_fig} for the 1-cover (left) and the 2-cover (right).}
				\label{pd_mc_fig}
\enf

Unfortunately, in the general setting of multiparameter persistence, a concise visual representation such as the one provided by the persistence diagram is still unknown in general. Nevertheless, a more algebraic encoding of this information can be achieved through the concept of  \emph{persistent Betti numbers}. More precisely, given a bifiltered simplicial complex $\K$, the \emph{persistent Betti numbers $\{\b_q^{\bb, \dd}\}_{\bb \le \dd}$} are defined as the rank invariants of the associated homology groups $\{H_q(K_{\bb})\}_{\bb \in \SSS}$, i.e., 
$$\b_q^{\bb, \dd} := \dim\big( \Im\big(H_q(K_{\bb}) \to  H_q(K_{\dd})\big)\big), $$
where we always work with $\Z/2$-coefficients. In a single-parameter setting, the persistent Betti number $\b_q^{b, d}$ has a clear connection with the persistence diagram: it counts the number of points that are contained in the upper-left domain $[0, b] \ti [d, \ff)$. In other words, this is the number of $q$-features that are born before time $b$ and live past time $d$. 
To simplify notation, we set henceforth
$\bqn := \bq(\K_{\cdot, n}).$

\section{Main results}
\label{res_sec}

The main conceptual contribution of this work consists in establishing the strong consistency and asymptotic normality for persistent Betti numbers in the setting of multi-parameter persistence. This will be achieved first in the \emph{scalar setting}, and then, under stronger conditions, in the \emph{functional setting}. While scalar results pertain to the asymptotic behavior of $\bqn$ for a \emph{fixed index pair} $\bb, \dd \in \SSS$, functional limit results allow for a description of $\bqn$ when seen as a stochastic process with varying $\bb$ and $\dd$.

Before stating the main results, we stress the importance of functional limit results. This is because already in the single-parameter setting, many of the most popular characteristics extracted from the persistence diagram cannot be expressed in terms of a single fixed persistent number. This applies in particular to the total persistence. In contrast, it was shown in \cite[Corollary 3.4]{svane} how a functional CLT in the single-parameter setting yields the asymptotic normality of the total persistence. Similarly, one immediate use case of the multi-parameter functional CLT is the asymptotic normality for the linear combination of total persistences taken with respect to different single-parameter filtrations that are derived from the bifiltration. We will revisit such examples in the discussion of the simulation experiments in Section \ref{sim_sec}.

Moreover, the functional setting is also the most relevant one in terms of the methodological contributions. While scalar strong consistency and asymptotic normality are essentially corollaries from known results in the single-parameter setting, substantial novel work is needed for the functional results.

%
%SCAL
%
\subsection{Scalar consistency and asymptotic normality}
\label{scal_sub}

As elucidated above, we first comment on the limiting behavior of $\bqn$ as $n \to \ff$ for fixed $\bb, \dd \in \mc \SSS$. More precisely, in Theorems \ref{lln_rank_thm} and \ref{clt_rank_thm} below, we extend the \emph{strong consistency} and \emph{asymptotic normality} derived for single-parameter persistent Betti numbers in \cite[Theorems 1.11 and 1.12]{shirai} to the multi-parameter setting. To ease notation, we henceforth write $X(B) := \#\{i \ge 1 \co X_i \in B\}$ for the number of points of $X$ contained in a set $B \su \R^d$.

%SCAL LLN
\bet[Strong consistency for the scalar rank invariant]
\label{lln_rank_thm}
Let $q \le d$, $\bb \le \dd$ and assume that $\{(X_i, M_i)\}_{i \ge 1}$ is ergodic and that $\E[X([0, 1]^d)^p] < \ff$ for all $p \ge 1$. Then, there exists a deterministic $\bbq$ such that 
$$\lim_{n \to \ff} \f1{n^d}\E[\bqn] = \bbq,$$
and, almost surely,
$$\lim_{n \to \ff}\f1{n^d}\bqn  = \bbq.$$
\ent

One of the strengths of Theorem \ref{lln_rank_thm} is that it holds for general ergodic point processes satisfying the moment condition. On the other hand, in the specific setting of binomial point processes where the points are sampled iid from a given density, one of the main results of \cite{lesnick} is a different LLN with respect to a highly refined homotopy interleaving distance. 

Next, we proceed with asymptotic normality.

%CLT
\bet[Asymptotic normality for the scalar rank invariant]
\label{clt_rank_thm}
Let $q\le d$, $\bb \le \dd$ and assume that $\{(X_i, M_i)\}_{i \ge 1}$ is an independently marked homogeneous Poisson point process. Then, as $n \to \ff$, 
$$\f{\bqn - \E\big[\bqn\big]}{n^{d/2}} \Rightarrow Z,$$
where $Z$ is a centered normal random variable.
\ent
We note that an incremental extension of the proof would establish asymptotic normality of linear combinations of the form $\sum_{j \le j_0}a_j \b_{q, n}^{\bb_j, \dd_j}$. Hence, the Cram\'er-Wold device \cite[Theorem 7.7]{pm} also gives the following corollary on the normality of the multivariate marginals.

\bet[Asymptotic multivariate normality for the scalar rank invariant]
\label{clt_rank_cor}
Let $j_0 \ge1$, $q\le d$, $\bb_j \le \dd_j$ for $j \le j_0$ and assume that $\{(X_i, M_i)\}_{i \ge 1}$ is an independently marked homogeneous Poisson point process. Then, as $n \to \ff$, 
$$\f{(\{\b_{q, n}^{\bb_j, \dd_j}\}_{j \le k}) - \E\big[\{\b_{q, n}^{\bb_j, \dd_j}\}_{j \le k}\big]}{n^{d/2}} \Rightarrow Z,$$
where $Z$ is a $k$-dimensional centered normal random vector.
\ent

To prove Theorems \ref{lln_rank_thm} and \ref{clt_rank_thm}, we note that for fixed $\bb \le \dd$, we can always interpret $\bqn$ within the standard framework of single-parameter persistence. Indeed, $\bqn$ can be seen as the standard persistent Betti number indexed over a set consisting of two elements, namely $\{\bb, \dd\}$. Hence, the proof of the scalar limit theorems will reduce to discussing to what extent the conditions stated in \cite{shirai} are valid for the \v C- and the multicover-bifiltrations.

%FUNC
\subsection{Functional consistency and asymptotic normality}
\label{func_sub}
While Theorems \ref{lln_rank_thm} and \ref{clt_rank_thm} describe the large-volume limit $\bqn$ for fixed $\bb, \dd \in \SSS$, we now view $\bqn$ as an $\SSS^2$-indexed stochastic process in the variables $\bb$ and $\dd$. Already in the setting of single-parameter persistence this viewpoint is essential for statistical applications since the total persistence and many other of the most commonly used characteristics extracted from persistence diagram are not linear combinations of the persistent Betti numbers. However, they can be expressed as a continuous functional given all persistent Betti numbers as input.

Hence, also in the multi-parameter setting, it is essential to move beyond scalar consistency and asymptotic normality. Furthermore, in comparison to the single-parameter case, new difficulties appear since there is no longer a clear representation of the persistent Betti numbers through the persistence diagram. Thus, we can no longer speak of features that are born in certain time intervals and die at later points in time.

%FUNC CONS
First, to extend Theorem \ref{lln_rank_thm} to the process level, we impose additional continuity constraints on the point process and on the mark distribution. In particular, independently marked Poisson point processes will satisfy all these conditions. First, for $q \ge 1$ the \emph{reduced $q$th factorial moment measure} $\a_q^!$ of a stationary point process is determined by the disintegration formula
$$\E\Big[\hspace{-.3cm}\sum_{\substack{ i_1,\dots, i_q \ge 1 \\\text{pw. distinct}}}\hspace{-.5cm}f(X_{i_1}, \dots, X_{i_q})\Big] = \la \int_{\R^d}\int_{\R^{d(q - 1)}} \hspace{-.8cm}f(x_0, x_0 +  y_1, \dots, x_0 + y_{q - 1})\a_q^!(\d y_1, \dots, \d y_{q - 1}) \d x_0$$
for any measurable $f\co \R^{dq} \to \off$, see \cite[Definition 8.6]{poisBook}. Additionally, we recall that the distribution $\P^0$ of a \emph{typical mark}  is determined by the property that 
$$\E^0[f(M)] = \f1\la \E\Big[\sum_{X_i \in Q_1} f(M_i)\Big]$$
for any measurable $f\co \off \to \off$, see \cite[Definition 9.3]{poisBook}. 

%FUNC LLN
\bet[Strong consistency for the functional rank invariant]
\label{lln_proc_thm}
Let $T\ge 0$, $q \le d$ and assume that $\{(X_i, M_i)\}_{i \ge 1}$ is ergodic and that $\E[X([0, 1]^d)^p] < \ff$ for all $p \ge 1$. Moreover, assume that $\a_{q + 2}^!$ are absolutely continuous, and that the distribution of the typical mark does not contain atoms. Then, almost surely, $\{\bqn\}_{(\bb, \dd)}$ converges to $\{\bbq\}_{\bb, \dd}$ in the sup-norm of processes on $\SSS^2$.
\ent

%FUNC NORM
Finally, we improve the CLT to a functional statement. This result is most naturally formulated in the Skorokhod space of multi-parameter c\`adl\`ag functions \cite{bickel}. To that end, we require additional regularity conditions on the mark distribution. Moreover, we assume that the point process $\{(X_i, M_i)\}_{i \ge 1}$ is in the subcritical regime of continuum percolation \cite{cPerc}. That is, with probability 1, there does not exist an infinite sequence $X_{i_1}, X_{i_2}, \dots$ of pairwise distinct points such that $|X_{i_j} - X_{i_{j + 1}}| \le T$ for all $j \ge 1$. 

%FCLT
\bet[Asymptotic normality for the functional rank invariant]
\label{clt_proc_thm}
Let $q\le d$ and $\{(X_i, M_i)\}_{i \ge 1}$ be an independently marked homogeneous Poisson point process. We assume that the distribution function of the typical mark is H\"older continuous with parameter $5/8$. Moreover, we assume that the point process $\{(X_i, M_i)\}_{i \ge 1}$ is in the subcritical regime of continuum percolation. Then, as $n \to \ff$, as a process
$$\f{\bqn - \E\big[\bqn\big]}{n^{d/2}}\Rightarrow \mc Z$$
where $Z$ is a centered Gaussian process.
\ent

\section{Simulation study}
\label{sim_sec}
In this section, we illustrate the asymptotic normality derived in Theorem \ref{clt_proc_thm} at the hand of simulated point patterns. This serves two purposes. First, we illustrate that the asymptotic normality is already clearly visible on bounded sampling windows. Second, the tests shed light on the power of goodness-of-fit tests that can be derived from the asymptotic normality.

To that end, we build on the simulation set-up that from \cite{krebs}. We will discuss the marked \v C-bifiltration and the multicover bifiltration separately in Sections \ref{cech_sub} and \ref{multi_sub}, respectively.

\subsection{Marked \v C-bifiltration}
\label{cech_sub}
As a null model, we take a homogeneous Poisson point process with intensity $\la = 0.2$ in a $10 \ti 10 \ti 10$ sampling window. By the marking theorem \cite[Theorem 5.6]{poisBook}, this can be also seen as a 2D homogeneous Poisson point process with intensity $\la = 2$ in a $10 \ti 10$ sampling window that is endowed with iid $\ms{Unif}([0, 10])$-marks. 

In order to reflect effects coming from the multiparameter persistence, we start from the general bifiltration $\K_{(r_1, r_2), n}$ from \eqref{bifi_eq} and extract three specific unifiltrations. First, fixing $r_1$ and varying $r_2$ corresponds to the sublevel-filtration with respect to the marks. Second, fixing $r_2$ and varying $r_1$ recovers the \v C-filtration. In addition to considering only these two specific axes, a popular general approach is to work with linear combinations in the sense that given $a, b \ge 0$ a simplex $\s$ is contained in the filtration at level $r$ if and only if $\s \in \Cech_{r_1}\big(\{X_i \in [0, n]^d \co M_i \le r_2\}\big)$ for some $r_1, r_2 \ge 0$ with $a r_1 + br_2 \le r$. To present specific examples, we look at the combinations $b = 1$ and $a\in \{5, 10, 20\}$. 

%
%TEST STAT
%
After having extracted several unifiltrations from the bifiltration in \eqref{bifi_eq}, we can now compute the associated persistence diagrams. This provides us with an ample choice of summary statistics that have been proposed in literature. One of the most basic ones is the \emph{total persistence} $\sum_i (D_i - B_i)$, i.e., the sum of all lifetimes when $\{(B_i, D_i)\}_i$ describe the birth- and death times in the persistence diagram. By our main result, this statistic is asymptotically normal in large sampling windows provided that the radii are sub-critical regime. Moreover, taking the total persistence of the mark filtration and of the \v C-filtration yields a bivariate random vector that is also asymptotically normal under the null hypothesis.

%
%FILT TEST
%
\subsubsection{Asymptotic normality}
\label{norm_mark_sec}
We now illustrate that the normal approximation becomes accurate in moderately-sized windows and extends beyond the sub-critical regime. To that end, we generate $n = 1,000$ realizations of the null model and recenter/rescale the total persistence by the empirical mean/standard deviation. To visualize the bivariate persistence vector, we first use the Cholesky decomposition of the estimated covariance matrix to decorrelate the coordinates and then compute the squared norm of the resulting vector. Under the null hypothesis, asymptotically this quantity follows a $\chi^2$-distribution. The resulting histograms in Figure \ref{high_mark_fig} already provide strong evidence for the approximate normality in moderately-sized windows. In general, this impression is reinforced by the Q-Q plots. However, a closer inspection reveals that the mark-filtration may have a small tendency towards more pronounced right and less pronounced left tails.

\bef[!h]
\begin{center}
	 \ig[width=.94\textwidth]{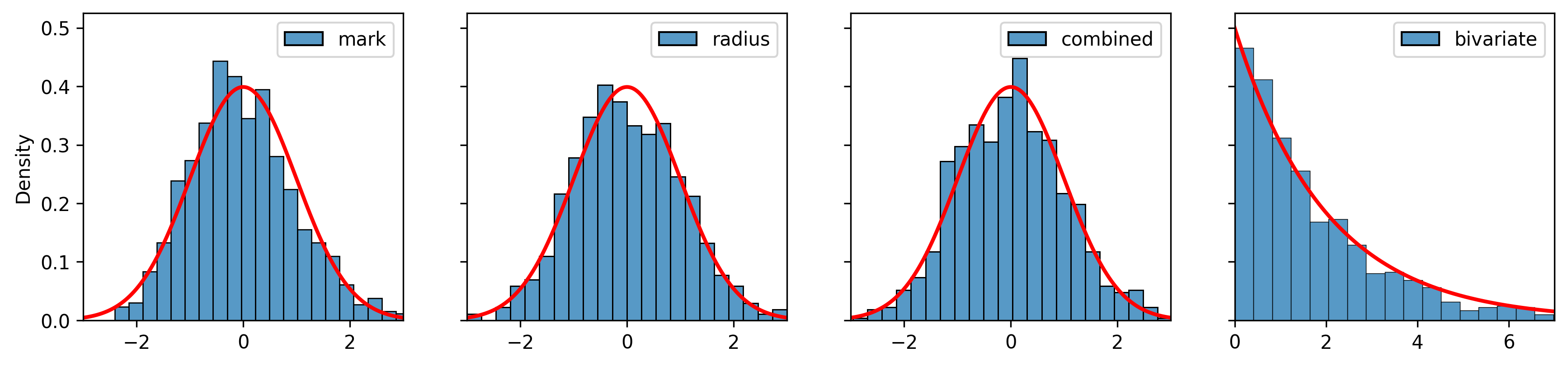}
	 \ig[width=.94\textwidth]{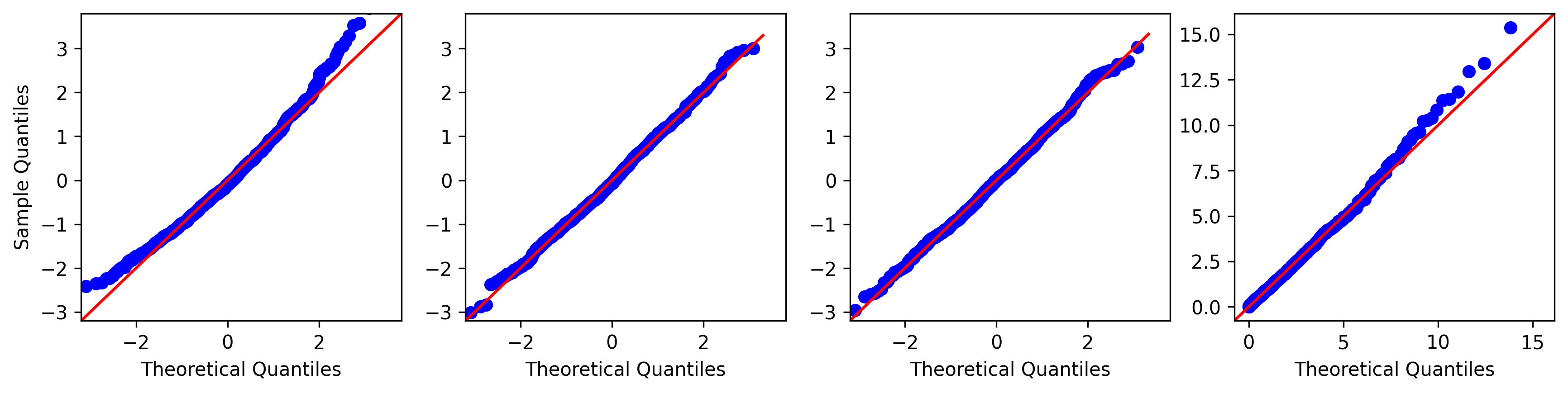}
	 \caption{Histograms and Q-Q plot (top/bottom row) for the total persistence in the mark-/\v C-/combined filtration (left/center left/center right, respectively). In the bivariate case, after decorrelation of the coordinates, the squared norm of the random vector is compared to a $\chi^2$ distribution.}
	 \label{high_mark_fig}
	 \end{center}
 \enf

 %
 %GOF TEST
 %
 \subsubsection{Goodness-of-fit tests}
 \label{good_mark_sec}
The most striking benefit of asymptotic normality are tests for the goodness of fit in large sampling windows. In this section, we illustrate this application at the hand of three specific point-process alternatives. They are prototypical examples for point patterns exhibiting attraction, repulsion, and a more complex interaction pattern, respectively.

For the attractive point pattern we work with a {Mat\'ern cluster process} based on a 3D parent Poisson process with intensity $0.2$, and where each parent generates a $\ms{Poi}(1)$-distributed number of offspring that are scattered uniformly  in a ball with radius 1. In particular, it has the same intensity as the null model. As repulsive pattern, we choose a {Strauss process} with intensity parameter $\b = 0.25$, interaction parameter $\g = 0.5$ and interaction radius $R = 1$. Finally, we also consider a {cell process} that exhibits a complex interaction structure but whose first- and second-order characteristics are similar to those of the null model. More precisely, we work with a modified Baddeley-Silverman process, where we partition the sampling window into $6^3$ congruent boxes and then place uniformly and independently in each of the boxes either 0, 1 or 2 points with probabilities 0.45, 0.1 and 0.45, respectively. Figure \ref{pattern_fig} shows realizations of the null model and the alternatives. It highlights that the differences between the models are subtle and difficult to recognize with the bare eye.

\bef[!h]
\centering	 \ig[width=1.09\textwidth]{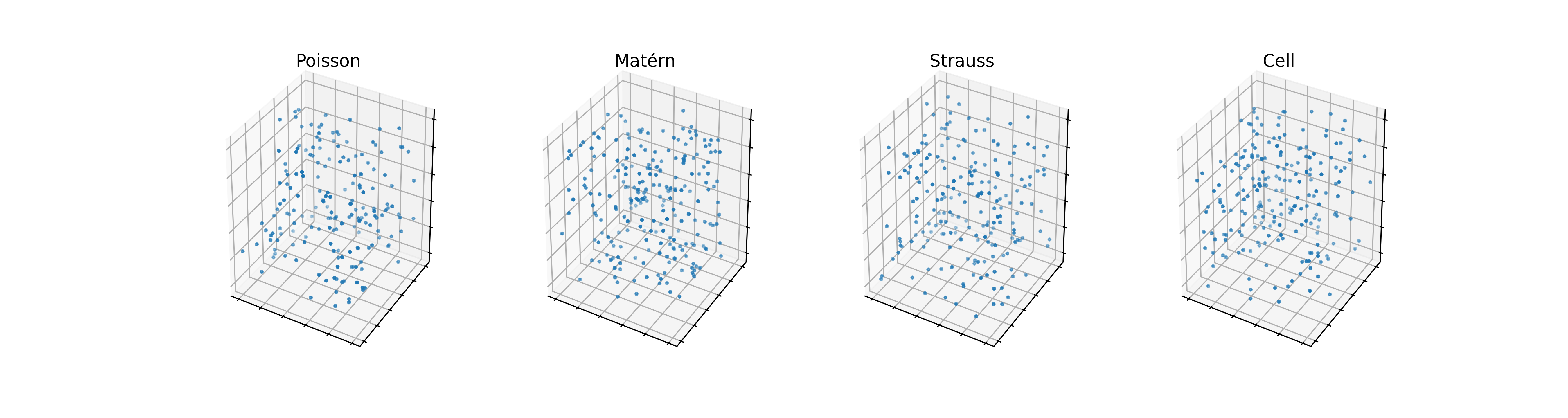}
	 \caption{Realizations of a Poisson, Mat\'ern cluster, Strauss, and cell process (from left to right).}
	 \label{pattern_fig}
 \enf

To perform the goodness-of-fit tests, we draw 1,000  samples from the Poisson null model in order to determine the mean and the standard deviation under the null model. Under the approximate normality, this yields the 5\% confidence region. Furthermore, for the bivariate setting, we proceed as in Section \ref{norm_mark_sec}. That is, under the null model we compute additionally a linear transformation that decorrelates the coordinates. Then, we compare the distribution of the squared vector norm with a $\chi^2$-distribution.

Based on also 1,000 samples from the alternatives, this makes it possible to compute the rejection rates, which are summarized in Table \ref{pow_mark_tab}. First, the type 1 error is close to the nominal level. 

Considering the rejection rates of the alternatives, we observe that the power of the test based on the combined filtration corresponding to $(a, b) = (20, 1)$ is superior to the pure mark, the pure \v C-filtration and also the corresponding bivariate vector. This illustrates one of the benefits of multiparameter persistence. 

We also observe that in general the test power is relatively low. On the one hand, this reflects that the alternatives are based on rather subtle deviations from the null hypothesis. However, a test based on Ripley's $K$-function illustrates that for the Mat\'ern and Strauss case, classical tools from spatial statistics are very powerful. Still when moving to the cell process, we see that such tests can become meaningless for point patterns involving complex interactions, whereas TDA-based tests still retain some power.

 \begin{table}[!htpb]
	 \begin{center}
		 \caption{Rejection rates for goodness-of-fit tests in the marked \v C-bifiltration.}
		 \label{pow_mark_tab}

		 \begin{tabular}{lllll}
			 & $\ms{Poi}$ & $\ms{Mat}$ & $\ms{Str}$ & $\ms{Cell}$ \\
			 \hline
 Total persistence -- mark filtration   & 4.4\% & 6.4\% & 15.2\% & 9.3\% \\
 Total persistence -- \v C filtration   & 5.2\% & 12.6\% & 13.7\%  & 15\%\\
			 Total persistence -- combined (5, 1)   & 5.4\% & 13.4\% & 15.7\%  &  14.3\%\\
			 Total persistence -- combined (10, 1)   & 5.8\% & 15.6\% & 17.1\%  &  14.9\%\\
			 Total persistence -- combined (20, 1)   & 5.7\% & 16.5\% & 17.4\%  &  15.2\%\\
			 Bivariate total persistence  & 4.6\% & 8.1\% & 14.3\%  &  12.4\%\\
 Ripley's $K$-function   & 5.1\% & 92.9\% & 92.8\%  &  3.4\%
		 \end{tabular}
	 \end{center}
 \end{table}

\subsection{Multicover bifiltration}
\label{multi_sub}

In contrast to Section \ref{cech_sub}, we set up the simulation study in two dimensions. The null model is a Poisson point process with intensity $\la = 2$ in a $10 \times 10$ sampling window.

As in Section \ref{cech_sub}, we extract two specific filtrations from the multicover filtration, namely the ones corresponding to the $k$-cover for $k \in \{1,2,3\}$. In particular, the case of the 1-cover reproduces the standard \v C-filtration on the underlying point pattern. As a test statistics, we also rely on total persistence. However, for the goodness-of-fit tests discussed in Section \ref{good_mc_sec} we found it to be beneficial to restrict to features born before a threshold. In order to provide a glimpse on the potential of characteristics combining the information inherent in different covers, we also consider a bivariate quantity where we form the weighted linear combination of the total persistences from the 1- and the 2-cover with weights $1/3$ and $2/3$, respectively. Again, by our main result, this statistic is asymptotically normal in large sampling windows provided that the radii are sub-critical regime. We now illustrate that the normal approximation becomes accurate in moderately-sized windows and extends beyond the sub-critical regime.

%
%FILT TEST
%
\subsubsection{Asymptotic normality}
\label{norm_mc_sec}

As in Section \ref{cech_sub}, we illustrate the asymptotic normality at the hand of histograms and Q-Q plots based on $n = 1,000$ realizations of the null model. Figure \ref{high_mc_fig} clearly indicates the asymptotic normality in the setting of the $k$-cover for any $k \in \{1, 2, 3\}$. 

\bef[!h]
\begin{center}
	 \ig[width=.94\textwidth]{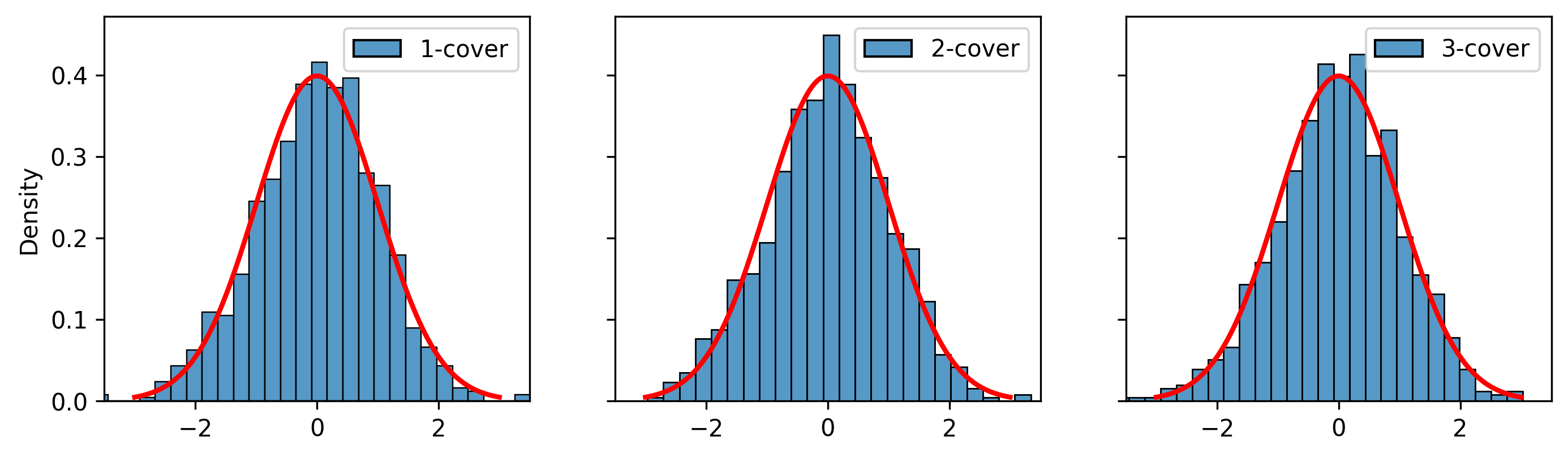}
	 \ig[width=.94\textwidth]{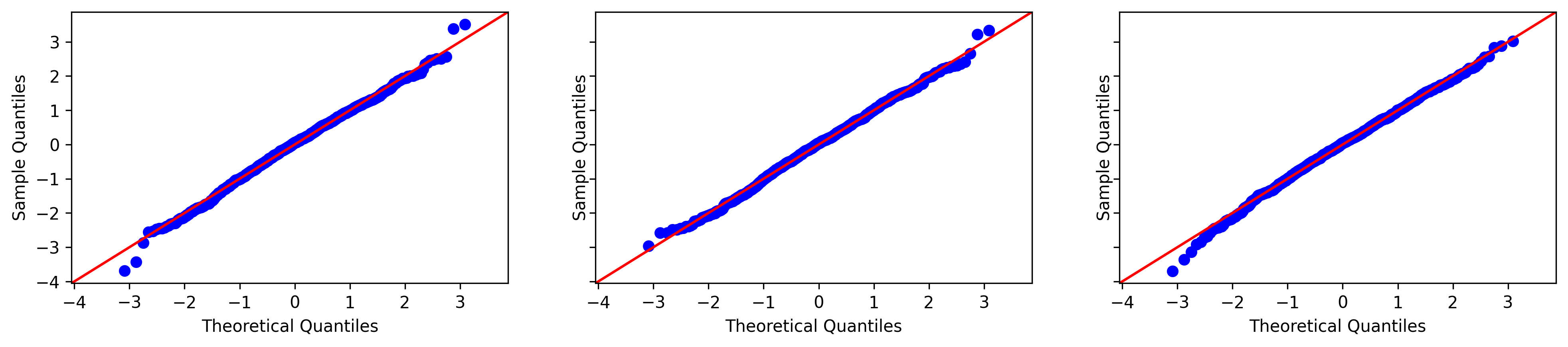}
	 \caption{Histograms (top row) and Q-Q plot (bottom row) for the total persistence in the \v C-filtration for the 1-cover (left) and 2-cover (right).}
	 \label{high_mc_fig}
	 \end{center}
 \enf

 %
 %GOF TEST
 %
 \subsubsection{Goodness-of-fit tests}
 \label{good_mc_sec}
 When evaluating the test power for possible alternatives, we work with examples that are similar to the ones considered in Section \ref{good_mark_sec} albeit we are now set in two dimensions.

 More precisely, as first alternative we fix a 2D {Mat\'ern cluster process} with parent-process intensity $0.2$ generating a $\ms{Poi}(1)$-distributed number of offspring uniformly scattered in a disk with radius 0.5. We also work again with a Strauss process, where we now choose $\b = 2.8$ as intensity parameter, $\g = 0.6$ as interaction parameter, and $R = 0.5$ as interaction radius. Finally, the cell process relies on a subdivision into 196 congruent boxes each containing either 0, 1 or 2 points with probabilities 0.45, 0.1 and 0.45, respectively.  Figure \ref{pattern_mc_fig} shows the patterns.

\bef[!h]
\hspace{-2cm}\ig[width=1.29\textwidth]{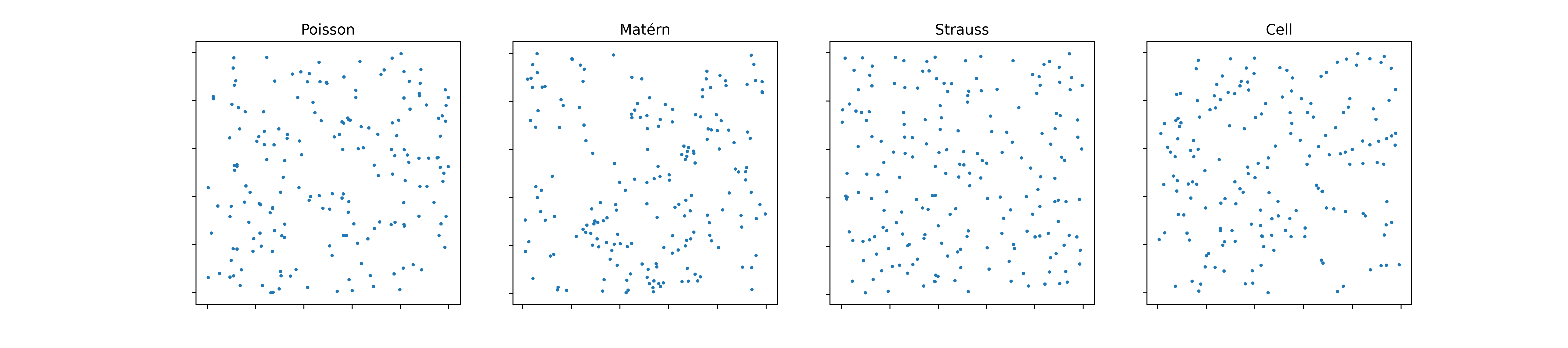}
	 \caption{Realizations of a Poisson, Mat\'ern cluster, Strauss, and cell process (from left to right).}
	 \label{pattern_mc_fig}
 \enf

We base the evaluation of the goodness-of-fit tests on 1,000 samples from the null model and from the alternatives, where we again rely on the asymptotic normality in order to construct the acceptance regions. First, the type 1 errors presented in Table \ref{pow_mc_tab} are close to the nominal 5\% level. When evaluating the test power at the alternatives, we see that the rejection rates are higher for the 2-cover filtration in comparison with the 1-cover filtration. Depending on the point pattern, the power for the test based on the 3-cover is sometimes higher and sometimes lower than that of the 2-cover. This gives again a glimpse on the potential of extracting different filtrations from multicover bifiltration. We also see that in general, the rejection rates are higher than the ones reported in Table \ref{pow_mark_tab}, although the latter is of course based on different point patterns, and in 3D.

Finally, a comparison with Ripley's $K$-function reinforces the impressions from Section \ref{good_mark_sec}: although classical characteristics from spatial statistics often exhibit superior performance for simple point patterns, they may fail to detect deviations from complete spatial randomness for point patterns with complex interactions. Here, TDA-based tests can offer valuable additional insights.

 \begin{table}[!htpb]
	 \begin{center}
		 \caption{Rejection rates for goodness-of-fit tests in the multicover bifiltration.}
		 \label{pow_mc_tab}

		 \begin{tabular}{lllcc}
			 & $\ms{Poi}$ & $\ms{Mat}$ & $\ms{Str}$ & $\ms{Cell}$ \\
			 \hline
 total persistence -- 1-cover   & 5.6\% & 38.3\% & 24.4\% & 13.3\% \\
 total persistence -- 2-cover   & 5.6\% & 63\% & 59.9\%  & 17.2\%\\
 total persistence -- 3-cover   & 4.7\% & 61.5\% & 63.4\%  & 8.7\%\\
total persistence -- Bivariate(1,2)  & 5.3\% &  65.4\% &  65.5\%  &  20.9\%\\
 Ripley's $K$-function   & 4.9\% & 93.4\% & 83.1\%  &  4.8\%
		 \end{tabular}
	 \end{center}
 \end{table}

\section{Conclusion and perspective}
\label{conc_sec}

Based on a functional LLN and a functional CLT for the rank invariant, the present work is the first step towards putting multi-parameter persistence on a sound statistical basis. However, in order to cover the variety of contexts where multi-parameter persistence is encountered in applications more research is needed. 

%MOD
First, on the modeling side, the assumption that the marks be independent of the location applies only rather specific application contexts. A more flexible framework would consider the setting of geostatistical markings, where the marks are determined through a possibly correlated random field in the background \cite{geost}. On a complementary side, it is also of interest to work with marks that are determined entirely by the underlying point configuration such as kernel-density estimators. The latter would present a refinement of the idea of the multicover bifiltration for delivering a persistence-based tool that is robust with respect to outliers.

%FUNC
Second, in multi-parameter persistence, the rank invariant is only a rather cumbersome summary statistics for the rich structure inherent in multi-parameter persistence modules. Recently, signed barcodes were proposed as a possibility to extend many of the strong points of the classical barcode representation in single-parameter persistence to the multi-parameter setting \cite{signed}. Hence, it would be worthwhile to investigate to what extent the asymptotic results of the present work extend to this characteristic.
Moreover, it would be exciting to develop statistics for the occurrence of specific indecomposable summands under different null hypotheses. The importance of this research direction for applications in materials science has also been stressed by Y.~Hiraoka in a recent series of talks and will be the topic of a manuscript \cite{shimizu}.

%POWER
Finally, the simulation in Section \ref{sim_sec} was designed as a proof of concept, where we illustrate that the asymptotic normality is already visible in bounded sampling windows, and where we provide first indications that  linear combinations of test statistics constructed from filtrations involving different parameters can outperform more classical single-parameter invariants. However, we have not touched upon the question on how to find a combination which is most powerful for discriminating the alternatives from the null model.  More generally, more research is needed in order to find out how to design test statistics that combine the information residing in the different layers of multi-parameter persistence in order to deliver the best performance for a given testing problem.

\section{Demonstration of asymptotic results}
\label{prel_sec}

Henceforth, we will often need to bound the change in the persistent Betti numbers $\bn$ when modifying $(\bb, \dd)$. In the following arguments, we rely on the representation
\begin{align*}
	\b_q^{\bb, \dd} &:= \dim\big( \Im\big(H_q(K_{\bb}) \to  H_q(K_{\dd})\big)\big) \\
	&= \dim(Z_q(K_{\bb})) - \dim(Z_q(K_{\bb}) \cap B_q(K_{\dd}))\\
	&= \dim(Z_q(K_{\bb}) + B_q(K_{\dd}))- \dim(B_q(K_{\bb}))  ,
\end{align*}
where $Z_q := \ker(\pa_q)$ and $B_q := \Im(\pa_{q + 1})$ denote the kernel and boundary spaces. Here, we tacitly identify $Z_q(K_{\bb})$ with its image in $Z_q(K_{\dd})$.  
We will need an extension of \cite[Lemma 2.11]{shirai}, which can be seen as a continuity result when changing the underlying filtration and filtration times.
We also write $K_{\bb}^q$ for the set of $q$-simplices in the simplicial complex $K_{\bb}$.

%
%CONT FILT LEM
%
\bel[Continuity of cycle and boundary spaces]
\label{cont_filt_lem}
Let $n \ge 1$ and $q \le d$. Let $\K = \{K_{\aa}\}_{\aa}$ and $\L = \{L_{\aa}\}_{\aa}$ be filtrations with $K_{\aa} \su L_{\aa}$ for all indices $\aa$. 
Moreover, let $\bb, \bb', \dd, \dd'$ be indices such that 
\been
\im $\bb\le \bb'$ and $\dd\le \dd'$;
\im $K_{\bb} \su K_{\bb'}$ and $K_{\dd} \su K_{\dd'}$.
\enen
Then, 
\been
\im $\dim(Z_q(L_{\bb'})) - \dim(Z_q(K_{\bb})) \le |L_{\bb'}^q \sm K_{\bb}^q|$; 
\im $\dim(B_q(L_{\bb'})) - \dim(B_q(K_{\bb})) \le |L_{\bb'}^{q + 1} \sm K_{\bb}^{q + 1}|$; 
\im $\dim(Z_q(L_{\bb'}) \cap B_q(L_{\dd'})) - \dim(Z_q(K_{\bb})\cap B_q(L_{\dd})) \le |L_{\bb'}^q \sm K_{\bb}^q| + |L_{\dd'}^{q + 1} \sm K_{\dd}^{q + 1}|. $
\enen
\enl
\bep
To ease notation, we set $Z' := Z_q(L_{\bb'})$, $Z:=Z_q(K_{\bb})$, $B' := B_q(L_{\dd'})$ and $B := B_q(K_{\dd})$.

For part 1., let $\psi$ and $\psi'$ be $q$-cycles in $Z'$ such that $\psi$ and $\psi'$ share the same $q$-simples in $L_{\bb'}^q \sm K_{\bb}^q$, then $\psi + \psi' \in Z$. Therefore,
$\dim(Z' ) - \dim(Z)\le  |L_{\bb'}^q \sm K_{\bb}^q|.$

For part 2., let $\pa(\psi), \pa(\psi') \in B'$ be $q$-boundaries coming from $(q + 1)$-chains $\psi, \psi'$ on $L_{\dd'}$ such that $\psi$ and $\psi'$ share the same simplices in $L_{\dd'}^{q + 1} \sm K_{\dd}^{q + 1}$. Then, $\psi + \psi'$ is a $(q + 1)$-chain in $K_{\dd}^{q + 1}$. Therefore,  
$\dim(B') - \dim(B) \le |L_{\dd'}^{q + 1} \sm K_{\dd}^{q + 1}|.$

The final claim now follows from the embeddings $(Z' \cap B') / (Z' \cap B) \su B' /B$ and $(Z' \cap B) / (Z \cap B) \su Z' /Z$.
\enp

For us, the main consequences of Lemma \ref{cont_filt_lem} are the continuity of the $\b$ and $\g$ characteristics.

%
%CONT FILT COR
%
\bec[Continuity of $\b$ and $\g$]
\label{cont_filt_cor}
Let $n \ge 1$ and $q \le d$. Let $\K = \{K_{\aa}\}_{\aa}$ and $\L = \{L_{\aa}\}_{\aa}$ be filtrations with $K_{\aa} \su L_{\aa}$ for all indices $\aa$. Moreover, let $\bb = (b_1, b_2, \kb), \bb'=(b_1', b_2', \kb), \dd = (d_1, d_2, \kd), \dd' = (d_1', d_2', \kd)$ be indices such that $\bb \le \dd$ and $\bb'\le \dd'$. Then, 
$\big|\bq(\K) - \bqp(\L)\big| \le  2|L_{\bb'}^q \sm K_{\bb}^q| +  |L_{\dd'}^{q + 1} \sm K_{\dd}^{q + 1}|.$
\enc

\section{Proof of Theorems \ref{lln_rank_thm} and \ref{clt_rank_thm}}
\label{scal_sec}

As announced in Section \ref{res_sec}, the main insight is that Theorems \ref{lln_rank_thm} and \ref{clt_rank_thm} may be interpreted in the setting of standard persistent homology by considering the binary filtration $\{\bb, \dd\} = \{(b_1, b_2, \kb), (d_1, d_2, \kd)\}$. We first reproduce the conditions (K1)--(K3) stated in \cite{shirai}. See also \cite{marked} for an extension to the marked setting. The assumption is that there exists a non-negative measurable function $\k$ defined on finite sets in $\R^d$ and taking values in $[0, \ff]$ such that:
\begin{itemize}
	\item[(K1)] $0 \le \kappa(\sigma) \le \kappa(\tau)$, if $\sigma$ is a subset of $\tau$;
	\item[(K2)] $\kappa$ is translation invariant, i.e., $\kappa(\sigma + x) = \kappa(\sigma)$ for any $x\in \R^d$;
		
	\item[(K3)] there is an increasing function 
		$\r \co [0, \ff) \to [0, \ff)$ such that
				       	\[|x - y| \le \r(\k(\{x,y\})),		\]
\end{itemize}
The $q$-simplices of the filtration at a level $r$ are then given by $(q + 1)$-subsets $\vp \su \R^d$ with $\k(\vp) \le r$.
Note that (K1)--(K3) are also extended to the marked settings in \cite{marked}.

In order to illustrate the idea in a simplified setting, we first consider the special case where we work with the 1-cover, i.e., where $k = 1$. Here, we first set $\k(\vp) := 3 + \diam(\vp)$ as soon as $\vp$ contains two elements that are distance at least $d_2$ apart. Otherwise, we distinguish on the number of elements of $\vp$. To compute $\bqn$ only the $q$-simplices and $(q + 1)$-simplices are relevant. Recall that $r_{\ms C}(\vp)$ denotes the filtration time in the \v Cech filtration.
If $\vp \su \R^d \times [0, T]$ consists of $q + 1$ elements, then we set $\k(\vp) := 0$ if $r_{\ms C}(\vp) \le b_1$ and all marks of the elements in $\vp$ are at most $b_2$. Otherwise, we set $\k(\vp) := 1$. Similarly, if $\vp \su \R^d \times [0, T]$ consists of $q + 2$ elements, then we set $\k(\vp) := 2$ if $r_{\ms C}(\vp) \le d_1$ and all marks of the elements in $\vp$ are at most $d_2$, and $\k(\vp) := 3$ otherwise. Finally, we set $\r(s) := d_2 + s$. Hence, the rank invariant $\bq$ corresponding to the birth time $\bb$ and death time $\dd$ in the multifiltration corresponds to the rank invariant for birth time $0$ and death time $2$ in the $\k$-filtration.

%
%MULTI-COV
%
For the general setting involving also the multicover filtration, the situation is a bit more delicate since the 0-simplices are no longer points in $\R^d$ but rather $k$-tuple of points. Hence, property (K3) does not make sense even syntactically any longer. However, the reason for introducing (K3) is to confine a simplex of filtration value $t > 0$ to a Euclidean ball with radius $\r(t)$, and this remains true in the multicover case: the diameter of any simplex at level $t$ is at most 2t.

\section{Proof of Theorem \ref{lln_proc_thm}}
\label{lln_sec}
%PROOF LLN PROC
To prove the functional convergence asserted in Theorem \ref{lln_proc_thm}, we build on a trick from \cite[Proposition 4.2]{thomas} and leverage the monotonicity properties of persistent Betti numbers. 

First, we explain how to derive Theorem \ref{lln_proc_thm} from a central continuity statement of the limiting persistent Betti number. 

On the geometric side, the key ingredient is the following continuity property of \v Cech filtration times from \cite[Lemma 6.10]{divol}.

\bel[Lipschitz continuity of \v Cech filtration times]
\label{divol_lem}
Let $ q \le d$ and $\rho > 0$. Let $X_0,\dots, X_q $ be iid uniform in a box $Q \su \R^d$.
Then, the function $F(r) := \P\big( r(X_0,\dots,X_q) \le r \big)$  is Lipschitz continuous on $[0, T]$.
\enl
Expressing the uniform distribution via the Lebesgue measure, we can also rewrite $F(r)$ in integral form
$$F(r) = |Q|^{-d(q + 1)}\int_{Q^{q + 1}} \one\{r(x_0, \dots, x_q) \le r\} \d (x_0, \dots, x_q).$$

In order to carry out monotonicity arguments, it will be convenient to extend the definition of $\bqn$ from indices $\bb$, $\dd$ satisfying $\bb \le \dd$ to general pairs by setting $\bqn := \b_{q, n}^{\bb, \max\{\bb, \dd\}}$, the maximum being taken pointwise.
%
%CONT LIM
%
\bel[Continuity of $\bar\b_q^{\bb, \dd}$]
\label{cont_lim_lem}
Assume that the distribution of the typical mark does not contain atoms. Then, $\bar\b_q^{\bb, \dd}$ is continuous in the indices $(\bb, \dd) \in \SSS^2$.
%, \bggq, \bgggq$ are 
\enl

Fixing a value $q \le d$, we suppress this parameter in the notation and simply write $\bar\b^{\bb, \dd}$  and $\b_n^{\bb, \dd}$.

%
%PRF LLN PROC
%
\bep[Proof of Theorem \ref{lln_proc_thm}]
Since the cover parameter $k \in \{1, \dots, T\}$ only takes a finite number of discrete values, we may fix $\kb, \kd \in \{1, \dots, T\}$ corresponding to the birth time and the death time for that parameter, respectively. Then, to simplify notation, for $\bb, \dd \in [0, T]^2$, we write $\b_n^{\bb, \dd}$ instead of $\b_n^{(\bb, \kb), (\dd, \kd)}$, and similarly $\bab^{\bb, \dd}$ instead of $\bab^{(\bb, \kb), (\dd, \kd)}$.
We split up the absolute value $|\b_n^{\bb, \dd} - \bab^{\bb, \dd}|$ into the positive and the negative part, which are considered separately.\\

\ni {$\bs{\sup_{\bb, \dd \in [0, T]^2} (\bab^{\bb, \dd } - \b_n^{\bb, \dd} )\to 0}$.}
Let $\e > 0$. By Lemma \ref{cont_lim_lem}, $\bab^{\bb, \dd}$ is continuous on the compact domain $[0, T]^4$ so that there exists $\de > 0$ such that $|\bab^{\bb, \dd} - \bab^{\bb', \dd'}| < \e$ whenever $|(\bb, \dd) - (\bb', \dd')|_\ff < \de$.

Note that $\b_n^{\bb, \dd}$ is increasing in $\bb$ and decreasing in $\dd$. Thus, writing $\SS_{\de, T} := [0, T]^2 \cap (\de\Z)^2$,
\begin{align*}
	        \sup_{\bb, \dd \in \SS^2} (\bab^{\bb, \dd}- \b_n^{\bb, \dd} )
			&=  \max_{\bb', \dd' \in\SS_{\de, T} }\sup_{(\bb, \dd) \in  (\bb', \dd') + [0, \de]^2\ti [-\de, 0]^2 }(\bab^{\bb, \dd} - \b_n^{\bb, \dd})\\
				&\le  \max_{\bb', \dd' \in\SS_{\de, T} }\sup_{(\bb, \dd)\in  (\bb', \dd') + [0, \de]^2\ti [-\de, 0]^2 }\Big((\bab^{\bb, \dd} - \bab^{\bb', \dd'}) + (\bab^{\bb' , \dd' } - \b_n^{\bb', \dd'})\Big)\\
	&, \e + \max_{\bb', \dd' \in\SS_{\de, T} }(\bab^{\bb' , \dd' } - \b_n^{\bb', \dd'}).
\end{align*}
Now, applying Theorem \ref{lln_rank_thm} shows that the maximum in the last line tends to 0 as $n \to \ff$, thereby concluding the proof.\\

\ni {$\bs{\sup_{\bb \le \dd \in [0, T]^2} (\b_n^{\bb, \dd} - \bab^{\bb, \dd })\to 0}$.} 
The argumentation is very similar to the setting considered above. More precisely, choosing $\bb', \dd' \in \de \Z^2$ such that $(\bb, \dd)$ is contained in $(\bb', \dd') + [-\de, 0]^2\ti [0, \de]^2$ gives that 
\begin{align*}
 \b_n^{\bb, \dd} -  \bab^{\bb, \dd}
						&\le  (\b_n^{\bb', \dd'} -  \bab^{\bb', \dd'}) +  (\bab^{\bb' , \dd' } - \bab^{\bb , \dd }) \le (\b_n^{\bb', \dd'} -  \bab^{\bb', \dd'}) + \e.
\end{align*}
We can now take $n$ large enough such that $(\b_n^{\bb', \dd'} -  \bab^{\bb', \dd'})$ is smaller than $\e$ uniformly over all $\de$-discretized indices.
\enp

The proof idea for Lemma \ref{cont_lim_lem} is to perform a reduction to the continuity of $n^{-d}\E[\b_n^{\bb, \dd}]$ in $(\bb, \dd)$ uniformly over all $n \ge 1$, which becomes a consequence of the stationarity.

%
%PROOF CONT LIM
%
\bep[Proof of Lemma \ref{cont_lim_lem}]
We may fix the discrete values $\kb, \kd \in \{1, \dots, T\}$ for the cover parameter corresponding to the birth time and the death time for that parameter, respectively.

By Theorem \ref{lln_rank_thm}, it suffices to show that $n^{-d}\E[\bqn]$ is continuous in each  $(\bb, \dd) \in [0, T]^4$ uniformly in $n \ge 1$. To that end, fix $\e > 0$ and $(\bb, \dd) \in [0, T]^4$. Now, in order to bound $\E[|\bn - \bnp|]$ for an arbitrary $(\bb', \dd') \in [0, T]^4$, we proceed stepwise and investigate the changes when modifying one of the indices at a time. 
More precisely, applying Corollary \ref{cont_filt_cor} twice,
$$n^{-d}\E[|\b_{q, n}^{\bb', \dd} - \bnp|] \le n^{-d}\big( \E\big[|K_{(\dd', \kd), n}^{q + 1} \De K_{(d_1', d_2, \kd), n}^{q + 1}|\big]+ \E\big[|K_{(\dd, \kd), n}^{q + 1} \De K_{(d_1', d_2, \kd), n}^{q + 1}|\big]\big),$$
where $\De$ denotes the symmetric difference. 
The main part of the proof is to show that there exists $\de > 0$ such that the right-hand side is smaller than $\e$ whenever $|d_1 - d_1'| \vee |d_2 - d_2'| \le \de$. This choice of $\de$ is uniform with respect to all $n \ge 1$ and all $d_1, d_1', d_2, d_2'$ satisfying $|d_1 - d_1'| \vee |d_2 - d_2'| \le \de$. Arguing in the same manner for the indices in $\bb'$ then concludes the proof.

%PALM
We start by explaining in detail how to proceed for $\kd = 1$. Afterwards, we elucidate how to argue for general $\kd \in \{1, \dots, T\}$. Moreover, we assume that $d_1 \le d_1'$, $d_2 \le d_2'$; the arguments are identical for the other possibilities.
Let $\P^0$ be the Palm probability for the marked point process $X = \{(X_i, M_i)\}_{i \ge 1}$, \cite[Definition 9.3]{poisBook}. That is, for any nonnegative measurable $f$,
$$\E^0[f(X)] = \f1\la\E\Big[\sum_{X_i \in [0, 1]^d} f(X - X_i)\Big].$$
If a $(q + 1)$-simplex $\s = (X_{i_0}, \dots, X_{i_{q + 1}})$ is contained in $K_{(\dd', \kd), n}^q \De K_{(\dd^\circ, \kd), n}^q$, then $r(\s) \in [d_1, d_1']$ or $M_{i_j} \in [d_2, d_2']$ for some $j \le q$. Writing $N_{1, n}, N_{2, n}$ for the number of ordered $(q + 1)$-simplices satisfying the first, respectively the second condition, it therefore suffices to show that there exists $\de > 0$ such that $n^{-d} \E[N_{i, n}] < \e$ whenever $d_i' - d_i < \de$. To achieve this goal, we note that
\begin{align*}
	\E[N_{1, n}] &\le \E\Big[\sum_{X_0 \in [0, n]^d}\#\{X_1, \dots, X_{q + 1} \in B_T(X_0) \text{ pw.~distinct}\co \rC(\s(X_0, \dots, X_{q + 1})) \in [d_1, d_1'] \}\Big]\\
	&\le \la n^d \int_{\R^{d(q + 1)}}\one\big\{\rC(\s(o, y_1, \dots, y_{q + 1})) \in [d_1, d_1']\big\} \a_{q + 1}^!(\d y_1,\d y_2, \dots, \d y_{q + 1}).
\end{align*}
Since $\a_{q + 1}^!$ is absolutely continuous with respect to the Lebesgue measure on $\R^{d(q + 1)}$, we deduce  from Lemma \ref{divol_lem} that there exists some $\de > 0$ such that $n^{-d}\E[N_{1,n}] < \e$ provided that $d_1' - d_1 < \de$. Similarly, any point $X_i \in X$ is contained in at most $X(B_T(X_i))^{q + 1}$ many $(q + 1)$-simplices of side length at most $T$. Hence, by the Cauchy-Schwarz inequality,
\begin{align*}
	\E[N_{2, n}] &\le  \E\Big[\sum_{X_0 \in [0, n]^d}X(B_T(X_0))^{q + 1}\one\{M_0 \in [d_2, d_2']\}\Big]\\
	&= \la n^d \E^0\big[X(B_T(o))^{q + 1} \one\{M_0 \in [d_2, d_2']\}\big]\\
	&\le \la n^d \sqrt{\E^0\big[X(B_T(o))^{2q + 2}\big]} \sqrt{\P^0(M_0 \in [d_2, d_2'])}
\end{align*}
so that again there exists some $\de > 0$ such that $n^{-d}\E[N_{2,n}] < \e$ provided that $d_2' - d_2 < \de$.

For general $\kd \in \{1, \dots, T\}$, we proceed essentially as above except that a simplex $\s$ is now described by $\kd (q + 1)$ not necessarily distinct points, i.e., 
$$\s =(X_{i_0, 1}, \dots,X_{i_0, \kd}, \dots, X_{i_{q + 1}, 1}, \dots,X_{i_{q + 1}, \kd}).$$
Still, the same computation as for the marked \v C-bifiltration shows that $n^{-d}\E[N_{i, n}]$ becomes small uniformly in $n$ provided that $d_i'$ is sufficiently close to $d_i$.
\enp

%
%FCLT
%
\section{Proof of Theorem \ref{clt_proc_thm}}
\label{fclt_sec}

The basis for the asymptotic normality of the persistent Betti numbers is a functional CLT. On a very general level, proving such a functional CLT involves two steps: 1) normality of the multivariate marginals, and 2) tightness. In Theorem \ref{clt_rank_thm} and the subsequent remark we established multivariate normality, so that it remains to verify tightness.

%NEW
To that end, we build on the blueprint used in \cite{krebs} and control moments of block increments. For multi-parameter persistence, three new challenges appear. First, \cite{krebs} was set in a quasi-1D domain, which simplified the cumulant expansion drastically. Second, more substantially, in multi-parameter persistence, we lose the interpretation of the increment $\b_n(E)$ as the number of features with birth- and death times in certain intervals because the algebraic structure of multiparameter persistence modules can be highly involved \cite{carlsson2009theory}. Third, the regularity conditions on the mark distribution need to be taken into account for in the moment bound.

%COV VS CECH
Although the general strategy applies to the combined model described in Example \ref{comb_exc}, writing out the proofs directly at this level of generality becomes cumbersome since 0-simplices in the $k$-cover filtration are represented by $k$-tuples of data points. Hence, we will generally first explain in detail the idea for $k = 1$, and then indicate how to adapt the arguments for the multicover bifiltration.

%TIGHT

In order to simplify notation, we henceforth fix the feature dimension $q \le d$ and do not highlight it further in the notation.
 A \emph{block} $E =\prod_{m \le 3}\Ebe \times \prod_{m \le 3}\Edd$ is a product of intervals $\Ebe,  \Edd \su [0, T]$. Writing $\Ebe = [b_{-, m}, b_{+, m}]$ and $\Edd = [d_{-, m}, d_{+, m}]$, we then define the increment of the rank invariant in the block $E$ as the alternating sum 
 \begin{align}
	 \label{alt_eq}
	 \b_n(E) := \sum_{i, i', i'', j,j',j'' \in \{-, +\}} {i i' i'' jj'j''}\b_n^{(b_{i, 1}, b_{i', 2}, b_{i'', 3}), (d_{j, 1}, d_{j', 2}, d_{j'', 3})},
 \end{align}
 where, $ii'i''jj'j''$ denotes the sign obtained when computing the product of the corresponding signs.
For instance, in the single-parameter case, i.e., if $\Ebe = \Edd = [0, T]$ for $m \in \{2, 3\}$, then $\b_n(E)$ describes the number of features that are born within the interval $[b_{-, 1}, b_{+, 1}]$ and die in the interval $[d_{-, 1}, d_{+, 1}]$. Note that there is an issue with the current definition of $\b_n(E)$ for instance, if $b_{+, 1} > d_{-, 1}$ since we would then need to compute an expression of the form $\b_n^{(b_{+, 1}, b_{-, 2}, b_{-, 3}), (d_{-,1}, d_{-,2}, d_{-,3})}$ which is ill-defined. However, we can actually exploit this situation to our advantage and set 
\begin{align*}
	\b_n^{(b_{+, 1}, b_{-, 2}, b_{-, 3}), (d_{-,1}, d_{-,2}, d_{-, 3})} &:= \b_n^{(b_{+, 1}, b_{-, 2}, b_{-, 3}), (b_{+,1}, d_{-,2}, d_{-, 3})} \\
	&\phantom{:=}+ \b_n^{(d_{-, 1}, b_{-, 2}, b_{-, 3}), (d_{-,1}, d_{-,2},d_{-,3})} \\
	&\phantom{:=}- \b_n^{(d_{-, 1}, b_{-, 2}, b_{-, 3}), (b_{+,1}, d_{-,2}, d_{-, 3})}.
\end{align*}
In particular, if $b_{-, 1} = d_{-, 1}$ and $b_{+, 1} = d_{+, 1}$, then $\b_n(E) = 0$. Hence, when deriving bounds on $\b_n(E)$, we may restrict to the cases, where $b_{+, m} \le d_{-, m}$. Note that this extension of $\b_n$ differs from the one considered in Section \ref{lln_sec}. However, since the claim in \ref{clt_proc_thm} only concerns the original and not the extended Betti number, this incompatibility does not cause any issues.

%CHENT
The Chentsov condition in the form of \cite[Display (3)]{bickel} shows that to verify tightness, it suffices to produce $C, \e >0$ such that for the centered increment $\bar\b_n(E) := \b_n(E) - \E[\b_n(E)]$ the moment condition
\begin{align}
	\label{lav_eq}
	n^{-2d}	\E[\bar\b_n(E)^4] \le C |E|^{1 + \e}
\end{align}
holds for all $n \ge 1$ and all blocks $E \su \SSS$. Note that if $E = E' \cup E''$ is the union of two adjacent blocks, then 
$$\bar\b_n(E)^4 = \bi(\bar\b_n(E') + \bar\b_n(E'')\bi)^4 \le 2^4 \max\{\bar\b_n(E')^4, \bar\b_n(E'')^4\}\le 2^4\bar\b_n(E')^4 + 2^4\bar\b_n(E'')^4.$$
For the discrete index corresponding to the multicover filtration, this observation implies that it suffices to treat the case where the intervals pertaining to that index reduce to singletons. That is, $\Ebbb = \{b_3\}$ and $\Edbb = \{d_3\}$.

 As in the proofs of \cite[Theorems 1, 2]{krebs}, we first show that it suffices to verify \eqref{lav_eq} for blocks belonging to a grid. In a second step, we derive crucial variance- and cumulant bounds for block increments in such a grid. Here, the \emph{mixed cumulant} of random variables $Z_1, \dots, Z_4$ with finite fourth moment is given by
\begin{align}
	\label{E:C4M1.1}
	        c^4(Z_1, \dots, Z_4)
		         = \sum_{\{N_1, \dots, N_r\}\prec\{1, \dots, 4\}}(-1)^{r - 1}(r - 1)!  \E\Big[\prod_{i \in N_1}Z_i\Big] \cdots \E\Big[\prod_{i \in N_r}Z_i \Big],
\end{align}
where the sum is taken over all the partitions $N_1 \cup \cdots \cup N_r$ of $\{1, \dots, 4\}$, see \cite[Identity (3.9)]{raic}. Moreover, we put $c^4(Z) := c^4(Z, Z, Z, Z)$.

In Sections \ref{grid_sec} and \ref{md_sec}, we prove the following key auxiliary steps. Henceforth, for a block $E = \prod_{m \le 3}\Ebe \times \prod_{m \le 3}\Edd$ we set $E_m := \Ebe \ti \Edd$. Then, we say that $E$ is \emph{$n$-good} if $|E_1| \ge n^{-d - \eb}$ and $|E_2| \ge n^{-8d/5 - \eb}$, where we set $\eb := 1/(100d)$. 

%
%GRID PROP
%
\bepr[Reduction to grid]
\label{grid_prop}
If the Chentsov condition \eqref{lav_eq} holds for all $n \ge 1$ and all $n$-good $E \su \oto$, then the processes $\{\bn\}_{\bb, \dd}$ are tight in the Skorokhod topology.
\enpr

Next, we state the cumulant bounds.

%
%VAR BOUND
%
\bepr[Variance and cumulant bound]
\label{var_prop}
It holds that 
$$\sup_{\substack{n \ge 1\\ E\su \oto \emph{ $n$-good}}}\f{\Var\big(\b_n(E)\big) + c^4\big(\b_n(E)\big)}{n^d |E_1|^{3/4 - \eb}|E_2|^{5/8 - \eb}}< \ff.$$
\enpr

We now explain, how to complete the proof of Theorem \ref{clt_proc_thm}. 
\bep[Proof of Theorem \ref{clt_proc_thm}]
First, we insert the bounds from Proposition \ref{var_prop} into the centered fourth moment expansion in terms of the variance and the fourth-order cumulant. Thus, for some $\ccv > 0$,
\begin{align*}
	\E[\bar\b_n(E)^4] &= 3\Var(\b_n(E))^2 + c^4(\b_n(E)) \\
	&\le 3\ccv^2 n^{2d} |E_1|^{3/2 - 2\eb}|E_2|^{5/4 - 2\eb} + \ccv n^d |E_1|^{3/4 - \eb}|E_2|^{5/8 - \eb}\\
	&\le 3\ccv^2 n^{2d}|E|^{5/4 - 2\eb} + \ccv|E|^{1 + \eb} n^d|E_1|^{-1/4 -2 \eb}|E_2|^{-3/8 -2 \eb}\\
	&\le 3\ccv^2 n^{2d}|E|^{5/4 - 2\eb} + \ccv|E|^{1 + \eb} n^dn^{d/4 + 4d\eb}n^{3d/5 + 4d\eb}.
\end{align*}
Since $n^dn^{d/4 + 4d\eb}n^{3d/5 + 4d\eb} \le n^{2d}$, we conclude the proof.
\enp

\subsection{Proof of Proposition \ref{grid_prop}}
\label{grid_sec}

%
%MOD CONT
%
By the analog of \cite[Theorem 16.8]{billingsley} in the multivariate setting, proving tightness requires a control on the modulus of continuity
$$\om'_\eta(\bar\b_n) := \inf_\La\max_{G\in\La} \sup_{(\bb, \dd), (\bb', \dd') \in L} |\bar\b_n^{\bb, \dd} - \bar\b_n^{\bb', \dd'}|,$$
where the infimum extends over all $\eta$-grids $\La$ in $\oto$.

To that end, we proceed in a similar vein as \cite[Proposition 5]{krebs}. More precisely, by invoking a result from \cite{davydov}, we show that it suffices to prove the claim for blocks from grids of the form $G_n = (n^{-\a_1}\Z \times n^{-\a_2}\Z)^2 $, $n \ge 1$ with $\a_1 := d/2 + \eb/2$ and $\a_2 := 4d/5 + \eb/2$. 
%
%GRID_PROOF MARKED CECH
%
\bep[Proof of Proposition \ref{grid_prop}]
To prove the claim, we will show that restricting $\bar\b_n^{\bb, \dd}$ to values $(\bb, \dd)$ taken from a sufficiently fine grid does not decrease the modified modulus of continuity substantially. More precisely, let $\e, \de >0$. Then, there exists $n_0 = n_0(\e, \de)$ such that almost surely
\begin{align}
	\label{dav_eq}
	\sup_{n \ge n_0}n^{-d/2}(\om'_\de(\bbn)- \om'_\de(\bbn|_{G_n})) \le \e.
\end{align}

%DAVYDOV
In order to prove \eqref{dav_eq}, we note that the process $\bn$ is increasing in $\bb$ and $(-\dd)$. Hence, by \cite[Corollary 2]{davydov}, it suffices to control
$$       \E\big[ \b_n^{(b_1', b_2, \kb), (d_1, d_2, \kd)} - \b^{(b_1, b_2, \kb), (d_1, d_2, \kd)}_n \big], \quad \E\big[\b_n^{(b_1, b_2', \kb), (d_1, d_2, \kd)} - \b^{(b_1, b_2, \kb), (d_1, d_2, \kd)}_n  \big] $$
for $|b_1' - b_1|\le n^{-\a_1}$, $|b_2' - b_2| \le n^{-\a_2}$ and 
$$\E\big[\b_n^{(b_1, b_2, \kb), (d_1, d_2, \kd)} - \b^{(b_1, b_2, \kb), (d_1, d_2', \kd)}_n \big], \quad \E\big[\b_n^{(b_1, b_2, \kb), (d_1, d_2, \kd)} - \b^{(b_1, b_2, \kb), (d_1', d_2, \kd)}_n \big]        $$
for $ | d_1 - d_1'|\le n^{-\a_1}$, $  | d_2 - d_2'| \le n^{-\a_2}$. By Corollary \ref{cont_filt_cor} it suffices to bound the number of simplices with filtration time contained in the corresponding interval of length $n^{-\a_i}$. We explain in detail how to proceed for $b_1' - b_1$ and $b_2' - b_2$, noting that the arguments for $d_1 - d_1'$ and $d_2 - d_2'$ are similar. Moreover, to make the arguments more accessible, we first lay out in detail how to proceed for the single-cover case, i.e., where $\kb = 1$. 

{$\boldsymbol{b_1' - b_1}.$} By applying the Mecke formula \cite[Theorem 4.4]{poisBook}, the expected number of $q$-simplices formed by Poisson points in $[0, n]^d$ is at most 
\begin{align*}
	&\E\big[\#\big\{X_{i_0}, \dots, X_{i_q} \in [0, n]^d \text{ pw.~distinct}\co r(X_{i_0}, \dots, X_{i_q}) \in [b_1, b_1']\big\}\big] \\
	&\quad \le \la^{q + 1}\int_{[0, n]^d} \int_{B_T(x_0)^q}\one\{r(x_0, x_1, \dots, x_q) \in [b_1, b_1']\} \d x_1\cdots \d x_q\d x_0\\
	&\quad \le \la^{q + 1}n^d \int_{B_T^q(o)}\one\big\{r(o, y_1, \dots, y_q) \in [b_1, b_1']\big\} \d y_1\cdots \d y_q.
\end{align*}
Hence, by Lemma \ref{divol_lem}, the expected number of $q$-simplices with \v Cech filtration time in $[b_1, b_1']$ is at most $$c\la^{q + 1}n^d|B_T(o)|^q (b_1' - b_1) \le c\la^{q + 1}|B_T(o)|^q n^{d - \a_1},$$
which is of order $o(n^{d/2})$.

{$\boldsymbol{b_2' - b_2}.$} Again, we apply the Mecke formula to bound the expected number of $q$-simplices with at least one vertex with mark in $[b_2, b_2']$ by 
\begin{align*}
	&\E\big[\#\big\{X_{i_0}\in [0, n]^d, X_{i_1}, \dots, X_{i_q} \in B_T(X_{i_0}) \text{ pw.~distinct}\co \{M_{i_0}, M_{i_1}, \dots, M_{i_q}\} \cap  [b_2, b_2'] \ne \es\big\}\big] \\
	&\quad \le \la^{q + 1}(q + 1)\int_{[0, n]^d} \int_{B_T(x_0)^q}\P(M \in [b_2, b_2']) \d x_1\cdots \d x_q\d x_0.
\end{align*}
Since the distribution function of the typical mark $M$ is H\"older-continuous, we deduce that the expected number of $q$-simplices with mark filtration time in $[b_2, b_2']$ is at most 
$$c\la^{q + 1}(q + 1) n^d|B_T(o)|^q (b_2' - b_2)^{5/8} \le c\la^{q + 1}(q + 1) n^d|B_T(o)|^q n^{d - 5\a_2/8},$$
which is of order $o(n^{d/2})$.

We now explain how to proceed when $\kb = k > 1$. \\

{$\boldsymbol{b_1' - b_1}.$} We need to bound the number of $q$-simplices  in the multicover filtration with filtration time in $[b_1, b_1']$. Any such $q$-simplex has the form 
$$\s =(X_{i_{0, 1}}, \dots, X_{i_{0, k}}, \dots, X_{i_{q, 0}}, \dots, X_{i_{q, k}}),$$
with $\rC(\s) \in [b_1, b_1']$ where the appearing points need not be distinct. Hence, an application of the Mecke formula gives the bound
$$\sum_{q \le \ell \le d} \E\big[\#\big\{X_{i_0}, \dots, X_{i_\ell} \in [0, n]^d \text{ pw.~distinct}\co \rC(X_{i_0}, \dots, X_{i_\ell}) \in [b_1, b_1']\big\}\big].$$
From here, we proceed as in the case $\kb = 1$.

{$\boldsymbol{b_2' - b_2}.$} Here, we need to bound the expected number of $q$-simplices with at least one vertex with mark in $[b_2, b_2']$. But again, any such $q$-simplex has the form
$$\s =(X_{i_{0, 1}}, \dots, X_{i_{0, k}}, \dots, X_{i_{q, 0}}, \dots, X_{i_{q, k}})$$
so that one of the associated marks $(M_{i_{0, 1}}, \dots, M_{i_{0, k}}, \dots, M_{i_{q, 0}}, \dots, M_{i_{q, k}})$ is contained in $[b_2, b_2']$. From here, we again proceed as in the case $\kb = 1$.
\enp

\subsection{Proof of Proposition \ref{var_prop}}
\label{md_sec}
In broad strokes, we follow the proof strategy of \cite[Theorems 1 \& 2]{krebs} and expand $\bar\b_n(E)$ as a sum of martingale differences. More precisely, we denote by $\{z_1, z_2, \dots, z_{\kn}\}$ the enumeration of $\{0, 1, \dots, n - 1\}^d$ according to the lexicographic order $\le_{\ms{lex}}$. Let
$$\GG_i := \s\big( \cup_{z \le_{\ms{lex}} z_i}(\PP \cap (z + [0, 1]^d))\big)$$
denote the $\s$-algebra generated by the Poisson points in cubes $z + [0,1]^d$ attached to lattice points $z \in \Z^d$ preceding $z_i$ in the lexicographic order. Then, $\b_n$ admits the martingale-difference decomposition
\begin{align}	\label{mdd_eq}
 \bar\b_n(E) = \sum_{i \le \kn}D_{i, n}(E),
\end{align}
where $D_{i, n}(E) := \E[\b_n(E)\ba \GG_i ] - \E[\b_n(E)\ba \GG_{i - 1}]$. In particular, $\Var(\b_n(E)) = \sum_{i \le \kn} \Var(D_{i, n}(E))$, so that to bound the variance,  it remains to control the second moments of $D_{i, n}(E)$ uniformly in $i \le \kn$ and $n \ge 1$ by a quantity of order $O(|E|^{1/2 + \e})$. Since the proof of Lemma \ref{diffBoundLem} is entirely analogous to \cite{krebs}, we omit the proof at this point.
%MOM LEM
\bel[Moment bound]        \label{diffBoundLem}
Let $k \ge 1$. Then,
$$\sup_{\substack{n \ge 1, i \le \kn\\ E\su \oto }}\f{\E[|D_{i, n}(E)|^k]}{ |E_1|^{3/4 - \eb/2}|E_2|^{5/8 - \eb/2}}< \ff.$$
\enl

%MOM
While the moment bounds from Lemma \ref{diffBoundLem} are enough to control the variance, dealing with the cumulants also requires the control of correlations. Since the proof of Lemma \ref{covBoundLem} can be copied from \cite{krebs}, we omit it at this point.
\bel[Covariance bound]
\label{covBoundLem}
For every $p_1, p_2 \ge 1$ there exist $C_{p_1, p_2} C_{p_1, p_2}' > 0$ with the following property. Let $n \ge 1$ and $A_1, A_2 \su \{1, \dots, \kn\}$  with $|A_1| = p_1$, $|A_2| = p_2$, and set          $X_{1, n}(E) = \prod_{i \in A_1}D_{i, n}(E)$ and $X_{2, n}(E) = \prod_{j \in A_2 }D_{j, n}(E)$. Then,
$$\Cov\big(X_{1, n}(E), X_{2, n}(E)\big) \le C_{p_1, p_2} \exp\big(-\dist(\{z_i\}_{i \in A_1},\{z_j\}_{j \in A_2})^{C_{p_1, p_2}'}\big)\sqrt{\E[X_{1, n}(E)^4]\E[X_{2, n}(E)^4]}.$$
\enl
We now elucidate how to complete the proof of Proposition \ref{var_prop} by extending the arguments from the cylindrical setup of \cite{krebs} to the general case considered here. After that, we prove Lemma \ref{diffBoundLem}.

\bep[Proposition \ref{var_prop}, cumulant bound]
We henceforth write $D_i$ for $D_{i, n}(E)$ and $\b_n$ for $\b_n(E)$.
Combining the martingale decomposition in \eqref{mdd_eq} with the multilinearity of cumulants yields that
\begin{align}\label{Cum_LargeBlocks}
	c^4(\bar\b_n(E)) = \sum_{ i ,j ,k ,\ell }a_{i, j, k, \ell} \ c^4\big(D_i, D_j, D_k, D_\ell\big),
\end{align}
where the $a_{i, j, k, \ell} \ge 1$ are suitable combinatorial coefficients, depending only on which of the indices $i, j, k, \ell$ are equal. In order to bound the right-hand side of \eqref{Cum_LargeBlocks}, we set $\de := \eb/1000$ and distinguish three cases. Every other case can be reduced to one of the three after permuting  $(i, j, k, \ell)$.
\been
\item\label{cum_a} $\diam(\{z_i, z_j, z_k, z_\ell\}) < |E|^{-\de}$,
\item\label{cum_b} $\dist(z_i, \{z_j, z_k, z_\ell\}) \ge |E|^{-\de}$,
\item\label{cum_c} $\dist(\{z_i, z_j\}, \{z_k, z_\ell\}) \ge |E|^{-\de}$ and $|z_i - z_j| \vee |z_k - z_\ell|\le |E|^{-\de}$.
\enen
\ni{\bf \eqref{cum_a}.} We need to bound the partial sum
$$\sum_{i,j,k,\ell \co \eqref{cum_a}}a_{i, j, k, \ell}  c^4(D_i, D_j, D_k, D_\ell).$$
Here, the sum is taken over all indices $i, j, k, \ell$ for which condition \eqref{cum_a} holds. To achieve this goal, we leverage the H\"older inequality and the representation in \eqref{E:C4M1.1}. Hence,
\begin{align*}
		        \big|c^4(D_i, D_j, D_k, D_\ell)\big| &\le \hspace{-.5cm}\sum_{\{N_1, \dots, N_r\}\prec\{i, \dots, \ell\}}\hspace{-.5cm}a_{\{N_1, \dots, N_r\}}'
				\prod_{m \in N_1}\hspace{-.1cm} \E[ |D_m|^{|N_1|} ]^{\f1{|N_1|}}  \cdots \prod_{m \in N_r}\hspace{-.1cm} \E[ |D_m|^{|N_r|} ]^{\f1{|N_r|}}\\
	&\le		        c \sup_{n \ge 1}\max_{\substack{k \le 4 \\ m \le n}}\E[|D_m|^k]
\end{align*}
with the coefficients $a_{\{N_1, \dots, N_r\}}'$ only depending on the structure of the partition. Now, by Lemma \ref{diffBoundLem} there exists $c' > 0$ such that $\E[|D_m|^k] \le c'|E_1|^{3/4 - \eb/2}|E_2|^{5/8 - \eb/2}$. Hence,
$$\sum_{\{N_1, \dots, N_r\}\prec\{i, \dots, \ell\}}a_{i, j, k, \ell}  c^4(D_i, D_j, D_k, D_\ell)\le 8cc'n|E|^{5/8 -\eb/2}|E|^{-3\de},$$
and the right-hand side is in $O(n|E|^{5/8 - \eb}).$

 \ni{\bf \eqref{cum_b} \& \eqref{cum_c}.} To avoid redundancy we provide only  the details for \eqref{cum_b}. That is, we control the sum
 $$\sum_{i,j,k,\ell \co \eqref{cum_b}}a_{i, j, k, \ell} c^4(D_i, D_j, D_k, D_\ell),$$
 To that end, we proceed along the lines of \cite[Proposition 9]{krebs}. More precisely, invoking the semi-cluster decomposition from \cite{barysh, raic} allows to express the individual cumulants in the form
 \begin{align*}
	c^4(D_i, D_j, D_k, D_\ell) =\hspace{-.5cm}\sum_{\{N_1, \dots, N_r\}\prec\{j, k, \ell\}} \hspace{-.5cm}a_{\{N_1, \dots, N_r\}}'\Cov\Big(D_i, \prod_{s \in N_1}D_s\Big)\E\Big[\prod_{s \in N_2}D_s\Big]\cdots \E\Big[\prod_{s \in N_r}D_s\Big]
 \end{align*}
 for some coefficients $a_{\{N_1, \dots, N_r\}}'$ only depending on the structure of the partition. Hence, combining the moment bounds from Lemma \ref{diffBoundLem} with the covariance bounds from Lemma \ref{covBoundLem} concludes the proof.
\enp

As in the previous proof, we henceforth write $D_i$ for $D_{i, n}(E)$ and $\b_n$ for $\b_n(E)$.

\bep[Proof of Lemma \ref{diffBoundLem}]
First, we rely on another representation of $D_{i, n}$. Namely, for $i \le n^d$, we set
$$\PP^{(i, n)} :=\big(\PP^{(n)}\sm (z_i +[0, 1]^d)\big) \cup \big(\wt{\PP}^{(n)}\cap  (z_i + [0,1]^d)\big),$$
where $\wt{\PP}^{(n)}$ is an independent copy of $\PP^{(n)}$. 
Then, $\E[\b_n\ba \GG_{i - 1}] = \E[\bnz \ba \GG_i]$, where $\bnz := \b_n(\PP^{(i, n)})$. Since, 
  $D_{i, n} = \E[\dnz \ba \GG_i ]$, where $\dnz :=  \b_n - \bnz$, it suffices to bound the moments of $|\dnz|$.
%
%CONN COMP 0
%
 First, we may restrict to simplices that are in a connected component intersecting $z_i + [0, 1]^d$. Indeed, let $\PP^{(n, *)}$ denote all points of $\PP^{(n)}$ whose connected component in $\bigcup_{X_i \in \PP^{(n)}}  B_T(X_i)$ intersects $z_i + [0, 1]^d$, and define $\PP^{(i, n, *)}$ similarly. Then,
$$\b_n = \b_n(\PP^{(n, *)}) + \bnz(\PP^{(n)} \sm \PP^{(n, *)}) \,\text{ and }\,\bnz = \bnz(\PP^{(i, n, *)}) + \bnz(\PP^{(i, n)} \sm \PP^{(i, n, *)}).$$
Since  $\PP^{(i, n)} \sm \PP^{(i, n, *)} = \PP^{(n)} \sm \PP^{(n, *)},$ we conclude that the difference can be computed with respect to $\PP^{(n, *)}$ and $\PP^{(i, n, *)}$. Setting $K := 500$, $p:= Kk(q + 1)$ and using that $\b_n(\PP^{(n, *)})$ and $\bnz(\PP^{(i, n, *)})$ are bounded above by the number of $q$-simplices in $\PP^{(n, *)}$ and $\PP^{(i, n, *)}$, respectively, the H\"older inequality gives that 
\begin{align*}
	\E[|\dnz|^k] &\le \E\Big[ \one\{\dnz \ne 0\}\Big((\#\PP^{(n, *)})^{q + 1}+ (\#\PP^{(i, n, *)})^{q + 1}\Big)^k\Big] \\
&\le 2^k\P(\bnz(\PP^{(n, *)})  + \bnz(\PP^{(i, n, *)}) \ne 0)^{1 - 1/K}\E\big[(\#\PP^{(n, *)})^p+ (\#\PP^{(i, n, *)})^p\big]^{1/K}\\
	&\le 2^{k + 1}\P\big(\bnz(\PP^{(n, *)}) \ne 0\big)^{1 - 1/K}\E\big[(\#\PP^{(n, *)})^p\big]^{1/K}.
\end{align*}
Since we assumed $X$ to be in the sub-critical regime of continuum percolation, the expected value in the last line is bounded above by a finite constant not depending on $i$ or $n$.

%DNZ BOUND
It remains to bound the probability $\P\big(\bnz(\PP^{(n, *)}) \ne 0\big)$. Now, we define
$$N_{i, n} := \min\big\{\ell \ge 1\co \PP^{(n, *)} \su Q_\ell(z_i)\big\}$$
as the size of the smallest box $Q_\ell(z_i) := z_i + [-\ell/2, \ell/2]^d$ centered at $z_i$ containing $\PP^{(n, *)}$.  Then,
\begin{align*}
\P\big(\bnz(\PP^{(n, *)}) \ne 0\big)&= \sum_{\ell \ge 1}\P(\bnz(\PP^{(n, *)}) \ne 0, N_{i, n} = \ell) \\
	&\le \sum_{\ell \ge 1}\P(N_{i, n} = \ell)^{1/K}\P(\bnz(\PP^{(n, *)}) \ne 0, N_{i, n} = \ell)^{1- 1/K}.
\end{align*}
Since $X$ is in the sub-critical regime of continuum percolation, we deduce that the tail probabilities $\sup_{n \ge 1}\sup_{i \le \kn}\P(N_{i, n}  = \ell)$ decay exponentially fast in $\ell$. Hence, it suffices to show that for some $c > 0$,
\begin{align}
	\label{dnzl_eq}
	\P(\bnz(\PP^{(n, *)}) \ne 0, N_{i, n} = \ell) \le  c \ell^{d(2q + 5)} |E_1|^{3/4}|E_2|^{5/8}.
\end{align}

%DEF F_i
 The key to \eqref{dnzl_eq} is to show that 
\begin{align} 	\label{dnzl_r_eq}
	\{\bnz(\PP^{(n, *)}) \ne 0, N_{i, n} = \ell\} \su F_{1, \ell} \cap F_{2, \ell},
\end{align}
where
\begin{align*}
	F_{1, \ell} := \big\{&\big(\rC(\s), r(\s')\big) \in E_1 \text{ for some $q$- and $(q + 1)$-simplices $\s$, $\s'$ contained in $Q_\ell(z_i)$}\big\}
\end{align*}
describes the event that there is at least one $q$-simplex with filtration time in $[b_{1, -}, b_{1, +}]$ and at least one $(q + 1)$-simplex with filtration time in $[d_{1, -}, d_{1, +}]$, and where
$$F_{2, \ell} :=  \big\{(M_k , M_j) \in E_2 \text{ for some $X_k, X_j \in Q_\ell(z_i)$} \big\}$$
denotes the event of finding marks in $E_2$.

%PRF INCL
Now, we prove \eqref{dnzl_r_eq}. First, if there is no $q$-simplex with filtration time in $[b_{1, -}, b_{1, +}]$, then in the alternating sum \eqref{alt_eq}, the contributions from $b_{1, -}$ and from $b_{1, +}$ cancel. Similarly, if there is no $(q + 1)$-simplex with filtration time in $[d_{1, -}, d_{1, +}]$, then the contributions from $d_{1, -}$ and $d_{1, +}$ cancel. Thus, $\{\bnz(\PP^{(n, *)}) \ne 0, N_{i, n} = \ell\} \su F_{1, \ell}$, and by similar arguments, $\{\bnz(\PP^{(n, *)}) \ne 0, N_{i, n} = \ell\} \su F_{2, \ell}$. 

%PROB BOUND
After having established \eqref{dnzl_r_eq} it remains to show that 
$\P(F_{1, \ell} \cap F_{2, \ell}) \le c \ell^{d(2q + 5)}|E_1|^{3/4}|E_2|^{5/8}.$
As in the previous proofs, we first deal with the single-cover setting, i.e., where $E_3 = \{(1, 1)\}$. Under the event $F_{1, \ell} \cap F_{2, \ell}$ there exist $2q + 5$ Poisson points $(X_{i_0}, M_{i_0})$, \dots, $(X_{i_q}, M_{i_q})$, $(X_{i_0'}, M_{i_0'})$, \dots,$(X_{i_{q + 1}'}, M_{i_{q + 1}'})$ and $(X_k, M_k)$, $(X_j, M_j)$ all contained in $Q_\ell(z_i)$. One slight nuisance is that several of these points may coincide, which needs to be taken into account before applying the Mecke formula. We present the detailed argumentation in one prototypical case, noting that the other ones can be handled similarly. More precisely, we assume that $i_0 = i_0' = k$, $i_1 = i_1' = j$ and $i_j = i_j'$ for $2 \le j \le q$, whereas apart from these, all other indices are pairwise distinct. Letting $F_\ell^*$ denote the event that the original event $F_{1, \ell} \cap F_{2, \ell}$ occurs due to this form of configuration, we deduce that
\begin{align*}
	\P(F_\ell^*)
	&\le \E\Big[\#\Big\{{X_{i_0}, \dots, X_{i_q}, X_{i_{q + 1}' }\in Q_\ell(z_i) \text{ pw.~distinct}}\co \\
	&\phantom{\le \E\Big[\#\Big\{}\big(r(X_{i_0}, \dots, X_{i_q}), r(X_{i_0}, \dots, X_{i_q}, X_{i_{q + 1}'})\big) \in E_1,\,		 (M_{i_0}, M_{i_1})\in E_2\Big\}\Big]\\
	&= \la^{q + 2} \int_{Q_\ell(z_i)^{q + 2}}\one\big\{\big(r(x_{0}, \dots, x_q)), r(x_0, \dots,  x_{{q + 1}})\big) \in E_1\big\}
\P\big((M_0, M_1)\in E_2\big)\d x_0 \cdots \d x_{q + 1}\\
	&\le c_{\ms{Lip}}^2\la^{q + 2}|E_2|^{5/8} \int_{Q_\ell(z_i)^{q + 2}}\one\Big\{\big(r(x_{0}, \dots, x_q)), r(x_{0}, \dots,  x_{{q + 1}})\big) \in E_1\Big\}\d x_0 \cdots \d x_{q + 1}.
\end{align*}
After scaling by $\ell^{(q +2 )d}$ and letting $X_0', \dots, X_{q + 1}'$ to be iid uniform in $Q_1(z_i)$, the final integral becomes the probability that $\big(r(X_0', \dots, X_q'), r(X_0', \dots, X_q', X_{{q + 2}}')\big) \in E_1$.
 Thus, by \cite[Proposition 6]{krebs}, we arrive at the asserted
$\P(F_\ell^*) \le c\la^{q + 2}\ell^{(q + 2)d}|E_1|^{3/4}|E_2|^{5/8}$
for some constant $c > 0$.

%GEN
Finally, we discuss how to extend the above proof to the general setting where $E_3 = \{(k, k')\}$ for some general $k \ge k' \ge 1$. Then, the $q$-simplex $\s$ and the $(q + 1)$-simplex $\s'$ have the form
$$\s =(X_{i_{0, 1}}, \dots, X_{i_{0, k}}, \dots, X_{i_{q, 0}}, \dots, X_{i_{q, k}}),$$
and 
$$\s' =(X_{i_{0, 1}'}, \dots, X_{i_{0, k'}'}, \dots, X_{i_{q + 1, 0}'}, \dots, X_{i_{q + 1, k'}'}).$$
Hence, we can bound $\P(F_{1, \ell} \cap F_{2, \ell})$ by the same arguments as in the single-cover setting. The only difference is that now there are more combinatorial possibilities concerning which of the points are pairwise distinct. Still, all of these cases can be handled through an application of a combination of the Mecke formula with \cite[Proposition 6]{krebs}.
\enp

\bibliography{lit}
\bibliographystyle{abbrv}

\end{document}